\documentclass[11pt,a4paper,openany]{article}

\usepackage{amsmath,amsfonts,ams}
\usepackage{latexsym}
\usepackage{amssymb}
\usepackage{color}
\usepackage{amscd}
\usepackage{amssymb}
\usepackage{amsthm}
\usepackage{graphicx}

\numberwithin{equation}{section}
\newtheorem{proposition}{Proposition}[section]
\newtheorem{definition}{Definition}[section]
\newtheorem{lemma}{Lemma}[section]
\newtheorem{theorem}{Theorem}[section]

\newtheorem{remark}{Remark}[section]

\begin{document}
\title{Scattering theory below energy for the cubic fourth-order Schr\"{o}dinger equation}
\author{{Changxing Miao$^1$, Haigen Wu$^{2}$, Junyong Zhang$^3$ }\\
          {\small $^1$ Institute of Applied Physics and Computational Mathematics}\\
         {\small P. O. Box 8009,\ Beijing,\ China,\ 100088}\\
         {\small $^2$ School of Mathematics and Information Science, Henan Polytechnic University  }\\
        {\small Jiaozuo,\ Henan Province,\ China,\ 454000 } \\
        {\small $^3$ Department of Mathematics, Beijing Institute of Technology  }\\
        {\small\ Beijing,\ China,\ 100081 } \\
         {\small (miao\_changxing@iapcm.ac.cn,  wuhaigen@gmail.com,  zhangjunyong111@sohu.com ) }
         \date{}
        }
\maketitle

\begin{abstract}
We investigate the global existence and scattering for the cubic
fourth-order Schr\"{o}dinger equation $iu_t+\Delta^2u+|u|^2u=0$ in
the low regularity space $H^s(\R^n)$ with $s<2$. We provide an
alternative approach to obtain a new interaction Morawetz  estimate
and extend the range of the dimension of the interactive estimate in
Pausader \cite{P08} by modifying a tensor product method appeared in
\cite{CGT}. We combine interaction Morawetz estimates, energy
increments for the I-method to prove the result.
\end{abstract}

\begin{center}
 \begin{minipage}{120mm}
 \end{minipage}
 \end{center}

\section{Introduction}
This paper is concerned with the global well-posedness and
scattering theory in $H^s(\R^n)$ of the Cauchy problem for the
following defocusing cubic fourth-order Schr\"{o}dinger equation
\begin{equation}\label{Eq}
\begin{cases}
i\partial_{t} u+\Delta^2u+|u|^2u=0~~~ \mathrm{in}~~~\R\times\R^n\\
u|_{t=0}=u_0,
\end{cases}
\end{equation}
with initial data $u_0\in H^s(\R^n)$ for some $0<s<2$ and $5\leq
n\leq7$. And we remark that this equation is corresponded to energy
critical case and mass critical case when $n=8$ and $n=4$,
respectively.

Fourth-order Schr\"{o}dinger equations have been introduced by
Karpman \cite{Kar96} and Karpman and Shagalov \cite{Kar00} to take
into account the role of small fourth-order dispersion terms in the
propagation of intense laser beams in a bulk medium with Kerr
nonlinearity. Such fourth-order Schr\"{o}dinger equations have been
studied from the mathematical viewpoint in Fibich, Ilan and
Papanicolaou \cite{FIP} who describe various properties of the
equaion in the subcritical regime, with part of their analysis
relying on very interesting numerical developments. Related
references given are by Ben-Artzi, Koch, and Saut \cite{BKS} who
gave sharp dispersive estimates for the biharmonic Schr\"odinger
operator which lead to the Strichartz estimates for the fourth-order
Schr\"odinger equation, see also \cite{MZ, P07, P08}.
 Nonlinear Schr\"odinger equations with third or fourth order anisotropic term
 have been discussed in Bocchel
\cite{Bou08},  for other special fourth order nonlinear
Schr\"{o}dinger equation, please refer to \cite{S06, HHW06,
HHW05}.
 We refer also to Pausader \cite{P07} where the energy critical case for
radially symmetrical initial data is discussed and Miao, Xu and Zhao
\cite{MXZ09, MXZ} simultaneously and independently obtained
scattering theory for the radially symmetrical initial data by using
argument  developed in  Killip and Visan \cite{KV}.  Miao and Zhang
\cite{MZ} showed the global well-posedness of the general high order
Schr\"odinger equation with defocusing nonlinearity. We also can
refer to Pausader \cite{P08} for the aim of finding a more completed
result on the cubic fourth-order Schr\"{o}dinger equation for
initial data $u_0\in H^2$ without radial assumption. However, very
little seems to be known about the existence and scattering theory
for fourth-order Schr\"odinger equation with large initial data in a
below energy space.

In this article, we prove the global well-posedness and scattering
theory of \eqref{Eq} in the lower regularity space $H^s(\R^n)$,
$s<2$, and those extend the global existence theory and scattering
result of Pausader \cite{P08} to the low regularity space. A main
ingredient in this paper is the new interaction Morawetz estimate
for the fourth-order Schr\"odinger equation. Our main result is the
following:
\begin{theorem}\label{the1}
The initial value problem \eqref{Eq} is globally-well-posed from
data $u_0\in H^s(\R^n)$ when $s>s_0$ and $5\leq n \leq7$. In
addition, there is a scattering for these solutions. Here
\begin{equation*}
s_0=~\begin{cases} \frac{16(n-4)}{7n-24}\qquad n=5,6\\
\frac{45}{23}\qquad\qquad n=7.
\end{cases}
\end{equation*}
\end{theorem}

The $I$-team \cite{C04T} introduced a new interaction Morawetz
potential for the nonlinear Schr\"{o}dinger equation
\begin{equation}\label{IMor}
M[u(t)]:=\int_{\R^3}|u(t,x)|^2\Big(\int_{\R^3}\Im[\bar
{u}(t,y)\nabla u(t,y)]\cdot \frac{x-y}{|x-y|}dy\Big)dx.
\end{equation}
This is a generalization of the classical Morawetz potential, which
has been studied in many literatures especially regarding on the
dispersive property of the Schr\"{o}dinger equations \cite{Bo99a,
Bo99b, G00, L78S, M68, N99}. The above equation \eqref{IMor}
generates a new space-time $L^{4}_{t,x}$ estimate for the defocusing
Schr\"{o}dinger equation with the general power nonlinearity.
Incorporating this with the almost conservation law, they showed
that the scattering of the equation and relaxed the low regularity
assumption given in the previous work \cite{C02T}. Two important
conserved quantities of equation \eqref{Eq} are the mass and the
energy. The mass is defined by
\begin{equation}\label{Mass}
M(u):=\frac12\int_{\R^n}|u(t)|^2d x,
\end{equation}
and the $H^2(\R^n)$ solutions satisfy the following energy
conservation
\begin{equation}\label{Energy}
E(u)(t):=\int_{\R^n}\frac12|\Delta u(t,x)|^2+\frac14|u(t,x)|^4d
x=E(u)(t_0).
\end{equation}
However the energy \eqref{Energy} of the $\dot H^s(\R^n)$( $s<2$)
solution can be infinite. The almost conservation law approach
allows us to monitor the energy of $Iu$ instead of a rough solution
$u$, where $I$ is a smoothing operator approximating to the identity
as passing to limit argument. In \cite{C02T}, this approach yields
that $\|u\|_{H^s(\R^3)}$ is bounded polynomially in time for the
cubic Schr\"{o}dinger equation, and the solution is globally
well-posed if $s>\frac56$. The regularity threshold is loosened to
$\frac45$ in \cite{C04T} due to the above mentioned new $L^4_{t,x}$
space-time estimate. All of these show that the interaction Morawetz
inequality plays an effective role in solving the low regularity
problem. In this paper, we will establish an interaction Morawetz
estimate and the almost conservation law for the defocusing
fourth-order Schr\"{o}dinger equation in the framework of $I$-method
to prove the global existence and scattering theory. We remark that
the interaction Morawetz estimate obtained by Pausader \cite {P08}
only for $n\geq 7$. Inspired by \cite{CGT}, we provide an
alternative approach to get an  interaction Morawetz inequality and
extend to Pausader's result to $n\geq 5$. However, there are some
differences and difficulties to derive the interactive estimate for
the fourth-order Schr\"odinger equation. It is well-known that the
Schr\"{o}dinger equation satisfies two local conservation laws. The
first one is  the local mass conservation $\partial_t
T_{00}+\partial_j T_{0j}=0$ and the other one is local momentum
conservation $\partial_t T_{j0}+\partial_k T_{jk}=0$ where
$T_{00}=\frac12|u|^2$ is the mass density and
$T_{0j}=T_{j0}=\Im(\bar u\partial_j u)$ is the momentum density and
the quantity
$$T_{jk}=2\Re(\partial_ju\bar{\partial_ku})+\delta_{jk}\big(-\frac12\Delta(|u|^2)+\frac{p-1}{p+1}|u|^{p+1}\big)$$
is the momentum current or stress tensor \cite{T06}. However, the
two above local conservation laws do not hold for the fourth-order
Schr\"{o}dinger equation and we shall utilize a modification of the
argument in \cite{CGT} to obtain the interactive estimate for the
fourth-order Schr\"odinger equation. As a direct consequence, we can
easily show that the solution to \eqref{Eq} is global and scatters
in the energy space $H^2(\R^n)$.

The paper is organized as follows. In Section 2, we recall the
Strichartz estimate for the fourth-order Schr\"odinger equation and
prove a local well-posedness of \eqref{Eq} in $H^s(\R^n)$ for
$s>\frac{n}2-2$ by the standard fixed point theorem. Section 3
provides an alternative approach to obtain the new interaction
Morawetz
 estimate. In Section 4, we prove the almost conservation
law for \eqref{Eq} by  frequency interaction strategy in
\cite{C04T}. In Section 5, the almost conservation law, the
interaction Morawetz inequality and a scaled bootstrap argument give
a uniform bound on $\|u(t)\|_{H^s(\R^n)}$ and the finiteness of
$\||\nabla|^{-\frac{n-5}4}u\|_{L^4_{t,x}}$. The scattering assertion
follows from the uniform bounds.

We conclude this introduction by setting some notations that will be
frequently used in this paper. If $X,Y$ are nonnegative quantities,
we sometimes use $X\lesssim Y$ or $X=O(Y)$ to denote the estimate
$X\leq CY$ for some $C$. Pairs of conjugate indices are written as
$p$ and $p'$ with $1\leq p\leq\infty$ and $1/p+1/{p'}=1$. We denote
$L^r=L^r(\R^n)$ to be the usual Lebesgue spaces. For $I\subset\R$,
we define the space-time space $L_t^qL_x^r$ by
$$\|u\|_{L_t^q(I;L^r_x)}:=\Big(\int_I\big(\int_{\R^n}|u(t,x)|^rdx\big)^\frac{q}rdt\Big)^\frac1q$$
with the usual modification when either $q$ or $r$ are infinity.
When there is no risk of confusion, we may shortened this norm to
$L^q_tL^r_x$ for readability, or to $L^r_{t,x}$ when $q=r$.

The Fourier transform on $\mathbb{R}^n$ is defined by
\begin{equation*}
\aligned \widehat{f}(\xi):= \big( 2\pi
\big)^{-\frac{n}{2}}\int_{\mathbb{R}^n}e^{- ix\cdot \xi}f(x)dx ,
\endaligned
\end{equation*}
giving rise to the fractional differentiation operators
$|\nabla|^{s}$ and $\langle\nabla\rangle^s$,  defined by
\begin{equation*}
\aligned
\widehat{|\nabla|^sf}(\xi):=|\xi|^s\hat{f}(\xi),~~\widehat{\langle\nabla\rangle^sf}(\xi):=\langle\xi\rangle^s\hat{f}(\xi),
\endaligned
\end{equation*} where $\langle\xi\rangle:=1+|\xi|$.
This helps us to define the homogeneous and inhomogeneous Sobolev
norms
\begin{equation*}
\big\|f\big\|_{\dot{H}^s_p(\R^n)}:= \big\|
|\nabla|^sf\big\|_{L^p_x(\R^n)},~~\big\|f\big\|_{{H}^s_p(\R^n)}:=
\big\| \langle\nabla\rangle^sf\big\|_{L^p_x(\R^n)}.
\end{equation*}

We will also need the Littlewood-Paley projection operators.
Specifically, let $\varphi(\xi)$ be a smooth bump function adapted
to the ball $|\xi|\leq 2$ which equals 1 on the ball $|\xi|\leq 1$.
For each dyadic number $N\in 2^{\mathbb{Z}}$, we define the
Littlewood-Paley operators
\begin{equation*}
\aligned \widehat{P_{\leq N}f}(\xi)& :=
\varphi\Big(\frac{\xi}{N}\Big)\widehat{f}(\xi), \\
\widehat{P_{> N}f}(\xi)& :=
\Big(1-\varphi\Big(\frac{\xi}{N}\Big)\Big)\widehat{f}(\xi), \\
\widehat{P_{N}f}(\xi)& :=
\Big(\varphi\Big(\frac{\xi}{N}\Big)-\varphi\Big(\frac{2\xi}{N}\Big)\Big)\widehat{f}(\xi).
\endaligned
\end{equation*}
Similarly we can define $P_{<N}$, $P_{\geq N}$, and $P_{M<\cdot\leq
N}=P_{\leq N}-P_{\leq M}$, whenever $M$ and $N$ are dyadic numbers.
We will frequently write $f_{\leq N}$ for $P_{\leq N}f$ and
similarly for the other operators.

The Littlewood-Paley operators commute with derivative operators,
the free propagator, and the conjugation operation. They are
self-adjoint and bounded on every $L^p_x$ and $\dot{H}^s_x$ space
for $1\leq p\leq \infty$ and $s\geq 0$, moreover, they also obey the
following
 Bernstein estimates
\begin{eqnarray*}\label{bernstein}
 \big\| P_{\geq N} f \big\|_{L^p} & \lesssim & N^{-s} \big\|
|\nabla|^{s}P_{\geq N} f \big\|_{L^p}, \\
\big\||\nabla|^s P_{\leq N} f \big\|_{L^p} & \lesssim  & N^{s}
\big\|
P_{\leq N} f \big\|_{L^p},  \\
\big\||\nabla|^{\pm s} P_{N} f \big\|_{L^p} & \thicksim & N^{\pm s}
\big\|
P_{N} f \big\|_{L^p},  \\
\big\| P_{\leq N} f \big\|_{L^q} & \lesssim &
N^{\frac{n}{p}-\frac{n}{q}} \big\|
P_{\leq N} f \big\|_{L^p},  \\
\big\| P_{ N} f \big\|_{L^q} & \lesssim &
N^{\frac{n}{p}-\frac{n}{q}} \big\|P_{ N} f \big\|_{L^p},
\end{eqnarray*}
where  $s\geq 0$ and $1\leq p\leq q \leq \infty$.

\section{The Strichartz Estimates and Local Well-Posedness}

The Strichartz estimates involve the following definitions:
\begin{definition}
A pair of Lebesgue space exponents $(q,r)$ are called Schr\"odinger
admissible for $\R^{n+1}$, or denote by $(q,r)\in \Lambda_0$ when
$q,r \geq 2, (q,r,n)\neq(2,\infty,2)$, and
\begin{equation}\label{Sadmissible}
\frac2q=n\Big(\frac12-\frac1r\Big).
\end{equation}
\end{definition}
\begin{definition} A pair of Lebesgue space exponents $(\gamma,\rho)$ are called
biharmonic admissible for $\R^{n+1}$,or denote by
$(\gamma,\rho)\in \Lambda_1$ when $\gamma,\rho \geq 2,
(\gamma,\rho,n)\neq(2,\infty,4)$, and
\begin{equation}\label{Biadmissible}
\frac4{\gamma}=n\Big(\frac12-\frac1{\rho}\Big).
\end{equation}
\end{definition}
\begin{proposition}[Strichartz estimates for
Fourth-order Schr\"odinger \cite{BKS, MZ,P07,
P08}]\label{Strichartz} Let $s\geq0$. Suppose that $u(t,x)$ is a
(weak) solution to the initial value problem
\begin{equation*}
\begin{cases}
(i\partial_t+\Delta^2)u(t,x)=F(t,x), \quad(t,x)\in [0,T]\times \R^n,\\
u(0)=u_0(x),
\end{cases}
\end{equation*}
for some data $u_0$ and $T>0$. Then we have the Strichartz estimate,
for $(q,r),(a,b) \in \Lambda_0$
\begin{equation}\label{Strichartz1}
\big\||\nabla|^su\big\|_{L^q([0,T];L^r)}\lesssim
\big\||\nabla|^{s-\frac2q}u_0\big\|_{L^2}+\big\||\nabla|^{s-\frac2q-\frac2a}F\big\|_{L^{a'}([0,T];L^{b'})},
\end{equation}
and for $(\gamma,\rho), (c,d)\in \Lambda_1$
\begin{equation}\label{Strichartz3}
\|u\|_{L^{\gamma}([0,T];L^{\rho})}\lesssim
\|u_0\|_{L^2}+\|F\|_{L^{c'}([0,T];L^{d'})}.
\end{equation}
\end{proposition}
As a consequence of the Strichartz estimate \eqref{Strichartz1} and
Sobolev's inequality, we have that
\begin{equation}\label{Strichartz2}
\|\Delta u\|_{L^q([0,T];L^r)}\lesssim \|\Delta u_0\|_{L^2}+\|\nabla
F\|_{L^{2}([0,T];L^{\frac{2n}{n+2}})},
\end{equation}
where $(q,r)$ is an any biharmonic admissible pair as in
\eqref{Biadmissible}.

The local existence theorem of \eqref{Eq} is as follows.
\begin{proposition}[Local Well-Posedness]\label{Local} Given any initial data $u_0\in H^s(\R^n)$ with $\frac n
2-2<s\leq2$, then there exists $T>0$ and a unique solution $u\in
C([0,T];H^s(\R^n))\cap L^{q_0}([0,T];\dot
H^{s+\frac2{q_0}}_{r_0}(R^n))\cap L^{q_1}([0,T];L^{r_1}(R^n))$ of
\eqref{Eq} for $(q_0,r_0)\in \Lambda_0$ and $(q_1,r_1)\in
\Lambda_1$.
\end{proposition}
\begin{proof}The proof is carried out by the standard fixed point
theorem together with the Strichartz estimate. For the sake of the
convenience and completeness, we merely sketch the proof for the
subcritical case with $n\geq 5$. To solve the equation \eqref{Eq} is
equivalent to solve the following integral equation
\begin{equation*}
u(t)=e^{it\Delta^2}u_0-i\int_0^te^{i(t-\tau)\Delta^2}|u|^2u(\tau)\mathrm{d}\tau.
\end{equation*}
Let $$X_T=C([0,T];H^s(\R^n))\cap X_0\cap X_1$$ where
\begin{equation*}
X_0=\bigcap_{(q_0,r_0)\in\Lambda_0}L^{q_0}([0,T];\dot
H^{s+\frac2{q_0}}_{r_0}(R^n))
\end{equation*}
and
\begin{equation*}
X_1=\bigcap_{(q_1, r_1)\in \Lambda_1}L^{q_1}([0,T];L^{r_1}(R^n)).
\end{equation*}
Let us define the Strichartz norm which is adapt to the Strichartz
estimate in Proposition \ref{Strichartz}
\begin{equation*}\label{Strichartz norm}
\|u\|_{X_T}:=\|u\|_{C([0,T];H^s(\R^n))}+\sup_{(q_0,r_0) \in
\Lambda_0}\||\nabla|^{s+\frac2{q_0}}u\|_{L^{q_0}_tL^{r_0}_x(I\times\R^n)}+\sup_{(q_1,r_1)
\in \Lambda_1}\|u\|_{L^{q_1}_tL^{r_1}_x(I\times\R^n)},
\end{equation*}
and a set
\begin{equation}\label{resolution space}
X:=\{u: \|u\|_{X_T}\leq 4C\|u_0\|_{H^s}\},
\end{equation}
and then we choose the space $(X, d)$ with metric
$d(u,v)=\|u-v\|_{L^{\frac {16}n}([0,T];L^4)}$ as a resolution space.
Then we claim that the solution map
$$A: u\mapsto e^{it\Delta^2}u_0-i\int_0^te^{i(t-\tau)\Delta^2}|u|^2u(\tau)\mathrm{d}\tau$$
is well defined for all $u\in X$ when $T$ is small enough.

Now we prove this claim. Actually, the Strichartz estimate yields
that
\begin{equation*}
\|A(u)\|_{X_T}\leq
2C\|u_0\|_{H^s}+C\big\||\nabla|^{s-1}(|u|^2u)\big\|_{L^{2}([0,T];L^{\frac{2n}{n+2}})}+C\big\||u|^2u\big\|_{L^1_t([0,T];L^2_x)}.
\end{equation*}
Thus the claim is reduce to prove the following nonlinear estimate
for some $\alpha>0$:
\begin{equation}\label{2.1}
\||u|^2u\|_{L^1_t([0,T];L^2_x)}+\big\||\nabla|^{s-1}(|u|^2u)\big\|_{L^{2}([0,T];L^{\frac{2n}{n+2}})}\leq
T^\alpha\|u\|^3_{X_T}.
\end{equation}
On one hand, when $\frac{n}2-2<s\leq\min\big\{\frac{n}2-1,2\big\}$,
we have
\begin{align}\nonumber
\big\||u|^2u\big\|_{L_t^1L_x^2}\lesssim&
T^{\frac12(s-\frac{n}2+2)}\|u\|_{L_t^2L_x^\frac{2n}{n-4}}\|u\|_{L_t^\frac8{n-2-2s}L_x^n}^2\\\nonumber
\lesssim&
T^{\frac12(s-\frac{n}2+2)}\|u\|_{L_t^2L_x^\frac{2n}{n-4}}\|u\|_{L_t^\frac8{n-2-2s}\dot
W^{s+\frac{n-2-2s}4,\frac{4n}{n+2+2s}}}^2\\\label{a2.1} \lesssim&
T^{\frac12(s-\frac{n}2+2)} \|u\|_{X_T}^3.
\end{align}
When $\frac{n}2-1<s<2$ (only happens in $n=5$),
\begin{equation}
\big\||u|^2u\big\|_{L_t^1L_x^2}\lesssim
T^\frac12\|u\|_{L_t^2L_x^\frac{2n}{n-4}}\|u\|_{L_t^\infty
L_x^5}^2\lesssim
T^\frac12\|u\|_{L_t^2L_x^\frac{2n}{n-4}}\|u\|_{L_t^\infty
H^s}^2\lesssim T^\frac12\|u\|_{X_T}^3.
\end{equation}
 On the other hand, under
the assumption $\frac n2-2<s<2$, it follows from fractional chain
rule, H\"older's inequality and Sobolev embedding in the case when
$\frac n2-2<s<\min\big\{\frac{n}2-1,2\big\}$ that
\begin{equation*}
\begin{split}
\big\||\nabla|^{s-1}(|u|^2u)\big\|_{L^{2}([0,T];L^{\frac{2n}{n+2}})}&\leq
T^{\frac{2s-n+4}4}\big\||\nabla|^{s-1}u\big\|_{L^\infty([0,T];L^{\frac{2n}{n-2}})}\|u\|^2_{L^{\frac{8}{n-2-2s}}([0,T];L^{n})}\\
&\leq T^{\frac{2s-n+4}4}\|u\|_{L^\infty([0,T]; \dot H^s)}
\big\||\nabla|^{s+\frac{n-2-2s}4}u\big\|^2_{L^{\frac{8}{n-2-2s}}([0,T];L^{\frac{4n}{n+2+2s}})}\\
&\leq T^{\frac{2s-n+4}4}\|u\|_{L^\infty([0,T]; \dot
H^s)}\|u\|_{X_T}^3,
\end{split}
\end{equation*}
and in the case when $\frac{n}2-1<s<2$ (only happens in $n=5$)  that
\begin{equation*}
\begin{split}
\big\||\nabla|^{s-1}(|u|^2u)\big\|_{L^{2}([0,T];L^{\frac{2n}{n+2}})}&\leq
T^\frac12\big\||\nabla|^{s-1}u\big\|_{L^\infty([0,T];L^\frac{2n}{n-2})}\|u\|_{L^\infty([0,T];L^5)}^2\\&\leq
T^\frac12\|u\|_{L^\infty([0,T]; H^s)}^3.
\end{split}
\end{equation*}
 Keeping in mind the norm of $X_T$,
then \eqref{2.1} follows from the above estimate and \eqref{a2.1}.
It can be similarly argued that $A$ is a contraction under the
metric $d(u,v)$. The existence and uniqueness assertion in $(X, d)$
follow from the fixed point theorem. Therefore, we conclude the
proof of this local well-posedness proposition.
\end{proof}
\section{The Interaction Morawetz Estimate in dimension $n\geq 5$.}
We adopt the convention that repeated indices are summed
throughout this section. Also, for $f,g$ two differentiable
functions, we define the mass and the momentum brackets by
\begin{equation*}
\{f,g\}_m=\Im(f\bar g) \quad \text{and}\quad
\{f,g\}_p=\Re(f\nabla\bar g-g\nabla\bar f).
\end{equation*}Given a smoothing real valued
function $a(x)$, we define the Morawetz action $M_a(t)$ by
\begin{equation}\label{MorA}
M_a(t)=2\int_{\R^n}\partial_j a(x)\Im\big(\bar u(x)\partial_j
u(x)\big)dx.
\end{equation}
\begin{proposition}[The Variation Rate of Morawetz Action]\label{MorAR}
If $u$ solves \eqref{Eq}, then the Morawetz action $M_a(t)$
satisfies the identity
\begin{equation*}\begin{aligned}
\partial_t M_a(t)=&2\int_{\R^n}\Big(2\partial_{jk}\Delta a\partial_ju\partial_k\bar u-\frac12\big(\Delta^3a\big)|u|^2-
4\partial_{jk}a\partial_{ik}u\partial_{ij}\bar
u\\
&\qquad \qquad+\Delta^2a|\nabla u|^2-\partial_j a\{|u|^2u,
u\}_p^j\Big)dx.\end{aligned}
\end{equation*}
\end{proposition}
\begin{proof}
Note that $\Im(z)=-\Re(iz)$, then it follows from the equation
\eqref{Eq} that
\begin{equation*}
\begin{split} \Im(\partial_t \bar u\partial_j u)=\Re(-i\partial_t\bar
u\partial_ju)=-\Re\big((\Delta^2\bar u+|u|^2\bar u)\partial_ju\big),
\end{split}
\end{equation*}
and
\begin{equation*}
\begin{split}
\Im(\bar u\partial_j \partial_t u)=\Re(-i\bar u\partial_j\partial_t
u)=\Re\big(\partial_j(\Delta^2u+|u|^2 u)\bar u\big).
\end{split}
\end{equation*}
Hence, a direction computation yields that
\begin{equation}\label{3.1}
\begin{split}
\partial_t M_a(t)=&2\int_{\R^n}\partial_j a\Re\big(\bar u\partial_j\Delta^2 u-\partial_j u\Delta^2\bar
u\big)dx -2\int_{\R^n}\partial_j a\{|u|^2u, u\}_p^jdx\\
=&-2\int_{\R^n}\Delta a\Re\big(\bar u\Delta^2
u\big)dx-4\int_{\R^n}\Re\big(\partial_j a\partial_j \bar
u\Delta^2u\big)dx -2\int_{\R^n}\partial_j a\{|u|^2u, u\}_p^jdx\\
:=&I_1+I_2-2\int_{\R^n}\partial_j a\{|u|^2u, u\}_p^jdx
\end{split}
\end{equation}
On one hand, we can see that from the integration by part
\begin{equation}\label{3.2}
\begin{split}
I_1=2\Re\int_{\R^n}\big(-\frac12\Delta^3 a |u|^2+\Delta^2 a |\nabla
u|^2+\partial_{jk}\Delta a\partial_j\bar u\partial_k u-\Delta
a|\partial_{jk}u|^2\big)dx.
\end{split}
\end{equation}
On the other hand, after a long length and careful computation, we
also have that
\begin{equation}\label{3.3}
\begin{split}
I_2=-4\Re\int_{\R^n}\big(\partial_{ijk}a\partial_j\bar
u\partial_{ik}u+2\partial_{jk}a\partial_{ij}\bar
u\partial_{ik}u-\frac12 \Delta a|\partial_{jk}u|^2\big)dx.
\end{split}
\end{equation}
Observe that
\begin{equation*}
\begin{split}
\Re\int_{\R^n}\partial_{ijk}a\partial_j\bar
u\partial_{ik}udx=-\Re\int_{\R^n}\partial_{jk}\Delta
a\partial_j\bar
u\partial_{k}udx-\Re\int_{\R^n}\partial_{ijk}a\partial_{ij}\bar
u\partial_{k}udx
\end{split}
\end{equation*}
and
\begin{equation*}
\begin{split}
\Re\int_{\R^n}\partial_{ijk}a\partial_{ij}\bar
u\partial_{k}udx=\Re\int_{\R^n}\partial_{ijk}a\partial_j\bar
u\partial_{ik}udx,
\end{split}
\end{equation*}
thus it follows that
\begin{equation}\label{3.4}
\begin{split}
\Re\int_{\R^n}\partial_{ijk}a\partial_j\bar
u\partial_{ik}udx=-\frac12\Re\int_{\R^n}\partial_{jk}\Delta
a\partial_j\bar u\partial_{k}udx
\end{split}
\end{equation}
Collecting \eqref{3.1}-\eqref{3.4}, it finally yields Proposition
\ref{MorAR}.
\end{proof}
We now derive a correlation estimate that is very useful in studying
the global well-posedness and the scattering properties of
fourth-order Shr\"odinger equations.
\begin{proposition}[Correlation Estimate for Fourth-Order Schr\"odinger Equation]\label{IMorEs}
If $u$ solves \eqref{Eq} on $[0,T]$, then we have the following
interactive Morawetz estimate that
\begin{equation}\label{IMorE}
\begin{split}
\|u\|^4_{M([0,T])}\lesssim \sup_{[0,T]}\|u(t)\|_{\dot
H^{\frac12}}^2\|u(t)\|_{L^2}^2.
\end{split}
\end{equation}
where
$\|u\|_{M([0,T])}:=\||\nabla|^{-\frac{n-5}4}u(x)\|_{L^4([0,T];L^4)}$.
\end{proposition}
\begin{proof}
Inspired by \cite{CGT}, we also introduce tensor product to derive a
correlation estimate for the fourth-order Schr\"odinger equation.
Let $u$ be solution to
\begin{equation*}
(i\partial_t+\Delta^2)u=F(u)
\end{equation*}
in $n$-spatial dimensions and $v$ be solution to
\begin{equation*}
(i\partial_t+\Delta^2)v=F(v)
\end{equation*}
in $m$-spatial dimensions. Define the tensor product $w:=(u\otimes
v)(t,z)$ for $z$ in $$\R^{n+m}=\{(x,y):x\in \R^n, y\in\R^m\}$$ by
the formula $$(u\otimes v)(t,z)=u(t,x)v(t,y).$$ One can check that
$w=u\otimes v$ solves the equation
\begin{equation}\label{EqT}
(i\partial_t+\Delta^2)w=F(u)\otimes v+F(v)\otimes u
\end{equation}
where $\Delta^2\triangleq\Delta_x^2+\Delta_y^2.$ Now we define the
Morawetz action $M_a^{\otimes_2}(t)$ corresponding to $w=u\otimes v$
by
\begin{equation}
\begin{split}
M_a^{\otimes_2}(t)&=2\int_{\R^n\otimes\R^m}\nabla
a(z)\cdot\Im\big(\overline{u\otimes v}(z)\nabla (u\otimes
v)(z)\big)dz\\
&=2\int_{\R^n\otimes\R^m}\nabla a(z)\cdot\Im\big(\bar{w}(z)\nabla
(w)(z)\big)dz,
\end{split}
\end{equation}
where $\nabla=(\nabla_x,\nabla_y)$. We now repeat the process of
proving Proposition \ref{MorAR} but more complicated to reach our
purpose. Also, it follows from the equation \eqref{EqT} that
\begin{equation*}
\begin{split} \Im(\partial_t \bar w\partial_j w)=\Re(-i\partial_t\bar
w\partial_jw)=-\Re\big((\Delta^2\bar w+|u|^2 {\bar u\bar v}+|v|^2
{\bar u\bar v})\partial_jw\big),\end{split}\end{equation*} and
\begin{equation*}
\begin{split}
\Im(\bar w\partial_j \partial_t w)=\Re(-i\bar w\partial_j\partial_t
w)=\Re\big(\partial_j(\Delta^2w+|u|^2 uv+|v|^2 uv)\bar w\big).
\end{split}
\end{equation*}
Moreover, we have that
\begin{align}\nonumber
&\partial_t M_a^{\otimes_2}(t)\\\nonumber
=&2\int_{\R^n\otimes\R^m}\partial_j a\Re\big(\bar
w\partial_j\Delta^2 w-\partial_j w\Delta^2\bar w\big)dz
-2\int_{\R^n\otimes\R^m}\partial_j a\{|u|^2uv+|v|^2uv,
w\}_p^jdz\\\nonumber =&-2\int_{\R^n\otimes\R^m}\Big[\Delta
a\Re\big(\bar w\Delta^2 w\big)+2\Re\big(\partial_j a\partial_j \bar
w\Delta^2w\big)\Big]dz -2\int_{\R^n\otimes\R^m}\partial_j
a\{|u|^2uv+|v|^2uv, w\}_p^jdz\\\label{3.6}
:=&II_1+II_2-2\int_{\R^n\otimes\R^m}\partial_j a\{|u|^2uv+|v|^2uv,
w\}_p^jdz.
\end{align}
A directional computation of expanding $\Delta^2 w$ in $II_1$ yields
that
\begin{equation}\label{3.7}
\begin{split}
II_1=&-2\int_{\R^n\otimes\R^m}\Delta a\Re\big(\bar u(x)\Delta_x^2
u(x)|v(y)|^2\big)dz\\
&-2\int_{\R^n\otimes\R^m}\Delta a\Re\big(\bar v(y)\Delta_y^2
v(y)|u(x)|^2\big)dz.
\end{split}
\end{equation}
In addition, for our purpose, we split the second term into several
pieces as follows:
\begin{equation*}
\begin{split}
II_2=&-4\Re\int_{\R^n\otimes\R^m}\partial_j a\big(\partial_j \bar
u(x)\bar v(y)+\bar u(x)\partial_j \bar
v(y)\big)\big(\Delta_x^2u(x)v(y)+u(x)\Delta_y^2v(y)\big)dz\\
=&-4\Re\int_{\R^n\otimes\R^m}\partial_j a\partial_j \bar
u(x)\Delta_x^2u(x)|v(y)|^2dz-4\Re\int_{\R^n\otimes\R^m}\partial_j
a\partial_j \bar u(x)\bar v(y)u(x)\Delta_y^2v(y)dz\\
&-4\Re\int_{\R^n\otimes\R^m}\partial_j a\bar u(x)\partial_j \bar
v(y)\Delta_x^2u(x)v(y)dz-4\Re\int_{\R^n\otimes\R^m}\partial_j
a\partial_j \bar v(y)\Delta_y^2v(y)|u(x)|^2dz\\
=&: II_2^{(1)}+II_2^{(2)}+II_2^{(3)}+II_2^{(4)}.
\end{split}
\end{equation*}
Observe that
\begin{equation*}
\begin{split}
II_2^{(2)}=2\Re\int_{\R^n\otimes\R^m}\Delta_x a\bar
v(y)\Delta_y^2v(y)|u(x)|^2dz
\end{split}
\end{equation*}
and
\begin{equation*}
\begin{split}
II_2^{(3)}=2\Re\int_{\R^n\otimes\R^m}\Delta_y a\bar
u(x)\Delta_x^2u(x)|v(y)|^2dz,
\end{split}
\end{equation*}
thus it follows that
\begin{equation*}
\begin{split}
&II_1+II_2^{(2)}+II_2^{(3)}\\
=&-2\int_{\R^n\otimes\R^m}\Big[\Delta_x a\Re\big(\bar u(x)\Delta_x^2
u(x)|v(y)|^2\big)+\Delta_y a\Re\big(\bar v(y)\Delta_y^2
v(y)|u(x)|^2\big)\Big]dz.
\end{split}
\end{equation*}
Hence, we can follow the computation of $I_1$ in \eqref{3.2} and see
that the right hand of the above can be written as
\begin{equation}\label{3.7}
\begin{split}
2\Re\int_{\R^n\otimes\R^m}\big(-\frac12\Delta_x^3 a |u|^2+\Delta^2_x
a |\nabla u|^2+\partial^x_{jk}\Delta_x a\partial_j\bar
u\partial_k u-\Delta_x a|\partial_{jk}u|^2\big)|v(y)|^2dz\\
+2\Re\int_{\R^n\otimes\R^m}\big(-\frac12\Delta_y^3 a
|v|^2+\Delta^2_y a |\nabla v|^2+\partial^y_{jk}\Delta_y
a\partial_j\bar v\partial_k v-\Delta_y
a|\partial_{jk}v|^2\big)|u(x)|^2dz,
\end{split}
\end{equation}
where $\partial^x_{jk}$ denote the the second order derivative with
respective to $x_j$ and $x_k$.   We also follow the computation of
$I_2$ in \eqref{3.3}  to obtain that
\begin{equation}\label{3.8}
\begin{split}
II_2^{(1)}=2\Re\int_{\R^n\otimes\R^m}\big(-4\partial_{jk}^x a
\partial_{ij}\bar u\partial_{ik}u+\partial^x_{jk}\Delta_x a\partial_j\bar
u\partial_k u+\Delta_x a|\partial_{jk}u|^2\big)|v(y)|^2dz
\end{split}
\end{equation}
and
\begin{equation}\label{3.9}
\begin{split}
II_2^{(4)}=2\Re\int_{\R^n\otimes\R^m}\big(-4\partial_{jk}^y a
\partial_{ij}\bar v\partial_{ik}v+\partial^y_{jk}\Delta_y a\partial_j\bar
v\partial_k v+\Delta_y a|\partial_{jk}v|^2\big)|u(x)|^2dz.
\end{split}
\end{equation}
Collecting \eqref{3.7}-\eqref{3.9}, we get that
\begin{equation}\label{3.10}
\begin{split}
II_1+II_2=2\Re\int_{\R^n\otimes\R^m}&\Big\{2\big(\partial^x_{jk}\Delta_x
a\partial_j\bar u\partial_k u|v|^2+\partial^y_{jk}\Delta_y
a\partial_j\bar v\partial_k
v|u|^2\big)\\-&\frac12(\Delta_x^3+\Delta_y^3) a
|uv|^2+\big(\Delta^2_xa|\nabla u|^2|v|^2+\Delta^2_ya |\nabla
v|^2|u|^2\big)\\-&4\big(\partial_{jk}^x a
\partial_{ij}\bar u\partial_{ik}u|v|^2+\partial_{jk}^y a
\partial_{ij}\bar v\partial_{ik}v|u|^2\big)\Big\}dz.
\end{split}
\end{equation}
Observe that if $a(z)=a(x,y)=|x-y|$, we have for $n\geq 5$
\begin{equation*}
\begin{split}
\Delta_x a=\Delta_y a=(n-1)|x-y|^{-1},\\
\Delta^2_x a=\Delta^2_y a=-(n-1)(n-3)|x-y|^{-3},\\
\end{split}
\end{equation*}
\begin{equation*}\Delta^3_x a=\Delta^3_y a=
\begin{cases}\quad
C\delta(x-y),\quad\quad\quad\quad\quad\quad \quad\qquad\quad n=5,\\
\quad3(n-1)(n-3)(n-5)|x-y|^{-5},\quad n\geq 6,
\end{cases}
\end{equation*}
and
\begin{equation*}
\begin{split}
\partial^x_{jk}
a=\partial^y_{jk}
a=|x-y|^{-1}\big(\delta_{jk}-\frac{(x-y)_j(x-y)_k}{|x-y|^2}\big),\\
\partial^x_{jk}\Delta_x
a=\partial^y_{jk}\Delta_y
a=-(n-1)|x-y|^{-3}\big(\delta_{jk}-\frac{3(x-y)_j(x-y)_k}{|x-y|^2}\big).
\end{split}
\end{equation*}
Now, for $e\in\R^n$ a vector, and $u$ a function, we define
\begin{equation*}
\nabla_eu=(e\cdot\nabla u)\frac{e}{|e|^2}\quad\text{and}\quad
\nabla_e^{\perp}u=\nabla u-\nabla_eu.
\end{equation*}
Therefore, for $e=x-y$, we can see that
\begin{equation*}
\begin{split}
2\partial^x_{jk}\Delta_x a\partial_j\bar u\partial_k u&=-2(n-1)|x-y|^{-3}|\nabla_e^{\perp}u|^2+4(n-1)|x-y|^{-3}|\nabla_eu|^2\\
\Delta^2_x a |\nabla u|^2&=-(n-1)(n-3)|x-y|^{-3}|\nabla u|^2\\
-4\partial_{jk}^x a
\partial_{ij}\bar u\partial_{ik}u&=-4|x-y|^{-1}\sum_{i}\big(|\nabla\partial_i u|^2-|\nabla_e\partial_i
u|^2\big)\leq\frac{-4(n-1)}{|x-y|^3}|\nabla_eu|^2
\end{split}
\end{equation*}
where we make use of the following estimate in the last inequality,
as shown in Levandosky and Strauss \cite{L00S} and \cite{P07}
\begin{equation*}
\sum_{i}\big(|\nabla\partial_i u|^2-|\nabla_e\partial_i
u|^2\big)\geq\frac{(n-1)}{|x-y|^2}|\nabla_eu|^2.
\end{equation*}
Making a similar argument for other terms, we finally control
$\partial_t M_a^{\otimes_2}$ as follows
\begin{equation*}
\begin{split}
\partial_t M_a^{\otimes_2}\leq2\Re\int_{\R^n\otimes\R^m}&\Big\{-\frac{2(n-1)}{|x-y|^{3}}\big(|\nabla_e^{\perp}u|^2|v|^2
+|\nabla_e^{\perp}v|^2|u|^2\big)\\&-\frac12(\Delta_x^3+\Delta_y^3) a
|uv|^2-\frac{(n-1)(n-3)}{|x-y|^{3}}\big(|\nabla
u|^2|v|^2+|u|^2|\nabla v|^2\big)\\&-2\partial_j a\{|u|^2uv+|v|^2uv,
w\}_p^j\Big\}dz.
\end{split}
\end{equation*}
Again through dropping some negative terms, we can dominate the
right hand part by
\begin{equation*}
\begin{split}
2\Re\int_{\R^n\otimes\R^m}\big(-\frac12(\Delta_x^3+\Delta_y^3) a
|uv|^2-2\partial_j a\{|u|^2uv+|v|^2uv, w\}_p^j\big)dz.
\end{split}
\end{equation*}
Hence, we get that
\begin{equation*}
\begin{split}
\int_{0}^T\int_{\R^n\otimes\R^m}\big((\Delta_x^3+\Delta_y^3) a
|uv|^2+4\partial_j a\{|u|^2uv+|v|^2uv,
w\}_p^j\big)dzdt\lesssim\sup_{[0,T]}|M_a^{\otimes_2}|,
\end{split}
\end{equation*}
that is
\begin{equation*}
\begin{split}
\int_{0}^T\int_{\R^n\otimes\R^m}\Big((\Delta_x^3+\Delta_y^3) a
|uv|^2+2\Delta_x a|u|^4|v|^2+2\Delta_y
a|v|^4|u|^2\Big)dzdt\lesssim\sup_{[0,T]}|M_a^{\otimes_2}|.
\end{split}
\end{equation*}
Choosing $u=v$, we get in the case that $n=5$
\begin{equation}\label{3.11}
\begin{split}
\int_{0}^T\int_{\R^n}|u(x,t)|^4dxdt\lesssim\sup_{[0,T]}|M_a^{\otimes_2}|
\end{split}
\end{equation}
and in the case that $n\geq 6$
\begin{equation}\label{3.12}
\begin{split}
\int_{0}^T\int_{\R^n\otimes\R^n}\frac{|u(x,t)|^2|u(y,t)|^2}{|x-y|^5}
dxdydt\lesssim\sup_{[0,T]}|M_a^{\otimes_2}|.
\end{split}
\end{equation}
However, we can write that
\begin{equation}\label{3.13}
\begin{split}
\int_{0}^T\int_{\R^n\otimes\R^n}\frac{|u(x,t)|^2|u(y,t)|^2}{|x-y|^5}
dxdydt=\int_{0}^T\int_{\R^n}|u(x,t)|^2(|u|^2\ast\frac{1}{|\cdot|^5})(x)
dxdt.
\end{split}
\end{equation}
 Now we define for $n\geq6$ the integral operator
\begin{equation*}
|\nabla|^{-(n-5)}f(x):=\int_{\R^n}\frac{f(y)}{|x-y|^5}dy
\end{equation*}
By applying Plancherel's Theorem to \eqref{3.13}, we obtain that
\begin{equation*}
\begin{split}
\int_{0}^T\int_{\R^n\otimes\R^n}\frac{|u(x,t)|^2|u(y,t)|^2}{|x-y|^5}
dxdydt=\int_{0}^T\int_{\R^n}\widehat{|u|^2}(\xi)|\xi|^{-(n-5)}\widehat{|u|^2}(\xi)
d\xi dt,
\end{split}
\end{equation*}
and the right hand also can be written as follows
\begin{equation*}
\begin{split}
\int_{0}^T\int_{\R^n}\big||\nabla|^{-\frac{n-5}2}\big(|u(x)|^2\big)\big|^2
dx dt.
\end{split}
\end{equation*}
For the sake of simplicity, we combine the two estimates
\eqref{3.11} and \eqref{3.12} pretending that $|\nabla|^0$ is
identity operator to get that for $n\geq5$
\begin{equation}\label{3.14}
\begin{split}
\big\||\nabla|^{-\frac{n-5}2}\big(|u(x)|^2\big)\big\|_{L^2([0,T];L^2)}^2\lesssim\sup_{[0,T]}|M_a^{\otimes_2}|.
\end{split}
\end{equation}
It can be shown by using Hardy's inequality (for details see
\cite{C04T})that for $n\geq5$
\begin{equation}\label{3.15}
\begin{split}
\sup_{[0,T]}|M_a^{\otimes_2}|\lesssim \sup_{[0,T]}\|u(t)\|_{\dot
H^{\frac12}}^2\|u(t)\|_{L^2}^2.
\end{split}
\end{equation}
 As so, it follows from \cite{P08,TVZ} that
$$\big\||\nabla|^{-\frac{n-5}4}u(x)\big\|_{L^4([0,T];L^4)}^4
\lesssim\big\||\nabla|^{-\frac{n-5}2}(|u(x)|^2)\big\|_{L^2([0,T];L^2)}^2.$$
This, together with \eqref{3.14} and \eqref{3.15}, yields that
\begin{align*}
\|u\|_{M[0,T]}^4=&\big\||\nabla|^{-\frac{n-5}4}u(x)\big\|_{L^4([0,T];L^4)}^4\\
\lesssim&\big\||\nabla|^{-\frac{n-5}2}(|u(x)|^2)\big\|_{L^2([0,T];L^2)}^2\\
\lesssim&\sup_{[0,T]}|M_a^{\otimes_2}|\\
\lesssim&\sup_{[0,T]}\|u(t)\|_{\dot H^\frac12}^2\|u(t)\|_{L^2}^2.
\end{align*}
\end{proof}

\begin{remark}
As an application of the interaction Morawetz estimate, one can
easily show that the solution to \eqref{Eq} is global and scatters
in the energy space $H^2(\R^n)$ with $5\le n\le 7$.
\end{remark}




\section{Almost Conservation Law}
The aim of this section is to control the growth in time of
$E(Iu)(t)$, where $Iu$ is a smoothing version of $u$. The operator
$I$ is a slightly modified smoothing operator as in \cite{C02T,
C04T}, depending on a parameter $N\gg1$ to be chosen later. For sake
of convenience, we recall the definition of operator $I$:
\begin{equation*}
\widehat{I f}(\xi):=m_N(\xi)\hat f(\xi),
\end{equation*}
where the multiplier $m_N(\xi)$ is smooth, radially symmetric,
nonincreasing in $|\xi|$ and
\begin{equation}\label{4.1}
m_N(\xi)=\begin{cases}1\qquad\qquad\qquad |\xi|\leq N,\\
\big(N|\xi|^{-1}\big)^{2-s}\quad |\xi|\geq 2N.\end{cases}
\end{equation}
Since $m_N(\xi)$ satisfies the H\"ormander multiplier condition, the
definition of $m_N(\xi)$ yields the following relationships between
$\|Iu\|_{H^2}$ and $\|u\|_{H^s}$ for $0<s<2$:
\begin{equation*}
\begin{split}
\|Iu\|^2_{H^2}\approx\int_{|\xi|\leq N}\langle \xi\rangle^{4}|\hat
u(\xi)|^2d\xi+\int_{|\xi|\geq 2N}\langle
\xi\rangle^{4}N^{2(2-s)}|\xi|^{-2(2-s)}|\hat u(\xi)|^2d\xi.
\end{split}
\end{equation*}
The right hand of the above can be controlled by
\begin{equation*}
N^{2(2-s)}\Big(\int_{|\xi|\leq N}\langle \xi\rangle^{2s}|\hat
u(\xi)|^2d\xi+\int_{|\xi|\geq 2N}|\xi|^{2s}|\hat
u(\xi)|^2d\xi\Big)\lesssim N^{2(2-s)}\|u\|^2_{H^s}.
\end{equation*}
Thus, we obtain that \begin{equation}\label{4.2}
\begin{split}
\|Iu\|^2_{H^2}\lesssim N^{2(2-s)}\|u\|^2_{H^s}.
\end{split}
\end{equation}
On the other hand, we can see that
\begin{equation*}
\begin{split}
\|u\|^2_{H^s}\lesssim\int_{|\xi|\leq N}\langle
\xi\rangle^{2s}|\widehat {Iu}(\xi)|^2d\xi+\int_{|\xi|\geq 2N}\langle
\xi\rangle^{2s}N^{-2(2-s)}|\xi|^{2(2-s)}|\widehat {Iu}(\xi)|^2d\xi.
\end{split}
\end{equation*}
The right hand of the above can be controlled by
\begin{equation*}
\int_{|\xi|\leq N}\langle \xi\rangle^{2s}|\widehat
{Iu}(\xi)|^2d\xi+\int_{|\xi|\geq 2N}|\xi|^{4}|\widehat
{Iu}(\xi)|^2d\xi\lesssim \|Iu\|^2_{H^2}.
\end{equation*}
Hence, the $L^2$ conservation yields that
\begin{equation}\label{4.3}
\|u\|^2_{H^s}\lesssim \|Iu\|^2_{\dot H^2}+\|Iu\|^2_{L^2}\leq
E(Iu)(t)+\|u_0\|^2_{L^2}.
\end{equation}
Once one has obtained a uniform bound on $E(Iu)(t)$ in terms of
$\|u_0\|_{H^s}$, the global well-posedness will follows from
\eqref{4.3}, the local well-posedness when $s>\frac n2-2$ and a
density argument. We remark a property of the operator $I$ in the
following lemma, which shows that the operator $\langle
\nabla\rangle I$ also holds the fractional Leibniz rule.
\begin{lemma}[ Leibniz rule]\label{Leibniz rule}
Let $1<r,r_1,r_2,q_1,q_2<\infty$ be such that
$\frac1r=\frac1{r_1}+\frac1{r_2}=\frac1{q_1}+\frac1{q_2}$ and the
$s$ in the operator $I$ satisfying $s\geq1$. Then
\begin{equation}
\big\|I\langle\nabla\rangle(fg)\big\|_{L^p}\lesssim\big\|(I\langle\nabla\rangle
f)\big\|_{L^{r_1}}\|g\|_{L^{r_2}}+\big\|(I\langle\nabla\rangle
g)\big\|_{L^{q_1}}\|f\|_{L^{q_2}}.
\end{equation}
\end{lemma}
The energy \eqref{Energy} is shown to be conserved
\begin{equation*}
\frac{d}{dt}E(u)(t)=\Re\int_{\R^n}\partial_t\bar
u(\Delta^2u+|u|^2u)dx=\Re\int_{\R^n}-i|u_t|^2dx=0.
\end{equation*}
Now we differentiate $E(Iu)(t)$ in time to obtain
\begin{equation*}
\begin{split} \frac{d}{dt}E(Iu)(t)&=\Re\int_{\R^n}\partial_t\overline{I
u}(\Delta^2Iu+|Iu|^2Iu)dx=\Re\int_{\R^n}\partial_t\overline{I
u}\big(|Iu|^2Iu-I(|u|^2u)\big)dx.
\end{split}
\end{equation*}
Integrating the above on the time interval $[0,t]$, we have
\begin{equation}\label{4.4}
\begin{split}E(Iu)(t)-E(Iu)(0)&=\Re\int_0^t\int_{\R^n}\partial_t\overline{I
u}\big[|Iu|^2Iu-I(|u|^2u)\big]dxdt:=\vartriangle E_I(t).
\end{split}
\end{equation}
Let us define $Z_I(t)$ as
\begin{equation}\label{4.5}
\begin{split}Z_I(t)=\sup_{(q,r) \in\Lambda_1}
\Big(\sum_{M\in 2^{\Z}}\|P_M\langle \Delta\rangle Iu\|^2_{L^q([0,t];L^r)}\Big)^{\frac12}.
\end{split}\end{equation}

Our aim is to show that the growth of $E(Iu)(t)$ satisfies
\begin{equation}
\vartriangle E_I(t)\lesssim N^{-\alpha}(Z_I(t))^\beta
\end{equation}
for some $\alpha,\beta>0.$  First, assuming a prior a small
space-time $\|u\|_M$ norm on the space-time slab $[0,t]\times \R^n$
and $E(Iu_0)$ is uniformly bounded, we can dominate $Z_I(t)$ in
terms of $\|u_0\|_{H^s}$.
\begin{lemma}\label{lem2}
Let $u(t,x)$ be in \eqref{Eq} defined on $[0,T^*]\times \R^n$ such
that
\begin{equation}\label{4.6}
\|u\|_{M([0,T^*])}\leq\delta
\end{equation}
for some small constant $\delta$. Assume $u_0\in C_0^\infty(\R^n)$
and $E(Iu_0)\lesssim1$. Then for $s\geq 1$ and $2>s>\frac n2-2$ and
sufficiently large $N$,
\begin{equation}\label{4.7}
Z_I(T^*)\leq C(\|u_0\|_{H^s(\R^n)}).
\end{equation}
\end{lemma}
\begin{proof}
Applying the operator $IP_M$, $\Delta IP_M$ to the equation
\eqref{Eq} and using the Strichartz estimate \eqref{Strichartz1} and
\eqref{Strichartz2},  we obtain for all $0\leq t\leq T^*$
\begin{align}\nonumber
Z_I(t)\lesssim&\|Iu_0\|_{H^2}+\big(\sum_{M\in2^{\mathbb{Z}}}\big\|P_M\nabla
I(|u|^2u)\big\|_{L_t^2L_x^\frac{2n}{n+2}}^2\big)^\frac12+\big(\sum_{M\in2^{\mathbb{Z}}}\big\|P_M
I(|u|^2u)\big\|_{L_t^2L_x^\frac{2n}{n+4}}^2\big)^\frac12\\\label{4.13}
\lesssim&\|Iu_0\|_{H^2}+\big\|\nabla
I(|u|^2u)\big\|_{L_t^2L_x^\frac{2n}{n+2}}+\big\|I(|u|^2u)\big\|_{L_t^2L_x^\frac{2n}{n+4}},
\end{align}
where we make use of Minkowski's inequality and the Littlewood-Paley
theory. By mean of the fractional chain rule, Lemma \ref{Leibniz
rule} and H\"older's inequality, we can control the nonlinearity as
follows
\begin{align}\label{4.14}
\big\|\nabla
I(|u|^2u)\big\|_{L_t^2([0,T^\ast];L_x^\frac{2n}{n+2})}\lesssim&
\big\|\nabla
Iu\big\|_{L^\infty([0,T^\ast];L_x^\frac{2n}{n-2})}\|u\|_{L^4([0,T^\ast];L_x^n)}^2,\\
\|I(|u|^2u)\|_{L_t^1([0,T^\ast];L_x^2)}\lesssim&\|Iu\|_{L_t^2([0,T^\ast];
L_x^\frac{2n}{n-4})}\|u\|_{L^4([0,T^\ast];L_x^n)}^2.
\end{align}
We write $$u=S_0 u+\sum_{j=1}^\infty\triangle_j u$$ where $S_0 u$
has spatial frequency support on $\langle \xi\rangle\lesssim N$ and
the remaining $\triangle_ju$ each have dyadic spatial frequency
support $\langle \xi\rangle \sim N_j:=2^{k_j}$, where $k_j\gtrsim
\log{N}$ is an integer for $j=1,2,\cdots$. Now we estimate
separately $\|u\|^2_{L^{4}([0,T^*];L^{n}_x)}$ on the low frequency
part $S_0 u$ and the high frequency pieces $\triangle_j u$, $j\geq
1$.

For the low frequency part $S_0 u$, the Sobolev embedding and
interpolation yield that
\begin{equation*}
\begin{split}
\|S_0 u\|_{L^{4}([0,T^*];L^{n}_x)}&=\|S_0
Iu\|_{L^{4}([0,T^*];L^{n}_x)}\leq \|Iu\|_{L^{4}([0,T^*];\dot
H^{\sigma}_{p})}\\&\leq \| Iu\|^{\theta}_{M([0,T^*])}\|
Iu\|^{1-\theta}_{L^{4}([0,T^*];\dot H^{2}_{\frac{2n}{n-2}})},
\end{split}
\end{equation*}
where
\begin{equation*}
\begin{split}
\theta&=\frac{2(8-n)}{7}>0,\qquad \sigma-\frac np=-1,\\
\sigma&=(\frac{5-n}4)\theta+2(1-\theta)=\frac{n^2-5n+4}{14}>0,\\
\frac1p&=\frac\theta4+\frac{n-2}{2n}(1-\theta)=\frac{n^2-5n+18}{14n}.\\
\end{split}
\end{equation*}
Thus, it follows that
\begin{equation}\label{4.15}
\|S_0 u\|_{L^{4}([0,T^*];L^{n}_x)}\lesssim
\delta^{\theta}{Z_I^{1-\theta}(T^*)}.
\end{equation}

For the high frequency pieces $\triangle_j u$ for $j=1,2,\cdots$,
the definition of $I$ operator gives that
\begin{equation*}
\|I\triangle_j u\|_{L^{4}([0,T^*];L^{n}_x)}\sim
N^{2-s}N_j^{s-2}\|\triangle_j u\|_{L^{4}([0,T^*];L^{n}_x)}.
\end{equation*}
Using the Bernstein inequality, we can rewrite that
\begin{equation*}
\|\triangle_j u\|_{L^{4}([0,T^*];L^{n}_x)}=N^{s-2}
N_j^{-s}N_j^{\frac{n-4}2}\|I\Delta\triangle_j
u\|_{L^{4}([0,T^*];L^{\frac{2n}{n-2}}_x)} .
\end{equation*}
Since $s>\frac n2-2$, we can sum in $j$ to obtain
\begin{equation}\label{4.16}
\sum_{N_j\geq N}\|\triangle_j u\|_{L^{4}([0,T^*];L^{n}_x)}\leq
N^{s-2}\sum_{N_j\geq N} N_j^{-[s-(\frac
n2-2)]}Z_I(T^*)=N^{\frac{n-8}2}Z_I(T^*).
\end{equation}
Collecting \eqref{4.13}-\eqref{4.16}, we conclude that
\begin{equation*}
\begin{split}
Z_I(t)\lesssim
\|Iu_0\|_{H^2}+\delta^{2\theta}{Z_I^{1+2(1-\theta)}(T^*)}+N^{{n-8}}Z^3_I(T^*).
\end{split}
\end{equation*}
Choosing $N$ sufficiently large and sufficiently small $\delta$, the
continuous argument yields \eqref{4.7} and thus it ends the proof of
Lemma \ref{lem2}.
\end{proof}

We have the following almost conservation law:
\begin{proposition}[Almost Conservation Law] \label{ACL}
Assume $2>s>\frac n2-2$ and $s\geq \frac{n-2}3$, $N\gg1$, $u_0\in
C_0^\infty(\R^n)$ and a solution of \eqref{Eq} on a time interval
$[0,T]$ for which
\begin{equation}\label{4.8}
\|u\|_{M([0,T])}\leq\delta
\end{equation} for some small constant $\delta$. In addition
we assume $E(Iu_0)\lesssim1$. Then we conclude that for all $t\in
[0,T],$
\begin{equation}\label{4.9}
E(Iu)(t)=E(Iu_0)+O(\max\{N^{\frac{n-8}2+},N^{-1+}\}).
\end{equation}
\end{proposition}
\begin{proof}
We apply the Parseval formula to $\vartriangle E_I(t)$ in
\eqref{4.4} to get
\begin{equation}\label{4.10}
\begin{split}
\vartriangle
E_I=\Re\int_0^T\int_{\sum_{j=1}^4\xi_j=0}\Big(1-\frac{m_N(\xi_2+\xi_3+\xi_4)}{m_N(\xi_2)m_N(\xi_3)m_N(\xi_4)}\Big)\\
\times \widehat{\overline{\partial_tI
u}}(\xi_1)\widehat{Iu}(\xi_2)\widehat{\overline{Iu}}(\xi_3)\widehat{Iu}(\xi_4)d\xi_2d\xi_3d\xi_4dt.
\end{split}
\end{equation}
Now if we use Equation \eqref{Eq} to substitute for $\partial_tI u$
in the above formula \eqref{4.10}, then it is split into two terms
as follows:
\begin{equation}\label{4.11}
\begin{split}
\vartriangle
E_1=\Big|\int_0^T\int_{\sum_{j=1}^4\xi_j=0}\Big(1-\frac{m_N(\xi_2+\xi_3+\xi_4)}{m_N(\xi_2)m_N(\xi_3)m_N(\xi_4)}\Big)\\
\times \widehat{\overline{\Delta^2I
u}}(\xi_1)\widehat{Iu}(\xi_2)\widehat{\overline{Iu}}(\xi_3)\widehat{Iu}(\xi_4)d\xi_2d\xi_3d\xi_4dt\Big|
\end{split}
\end{equation}
and
\begin{equation}\label{4.12}
\begin{split}
\vartriangle
E_2=\Big|\int_0^T\int_{\sum_{j=1}^4\xi_j=0}\Big(1-\frac{m_N(\xi_2+\xi_3+\xi_4)}{m_N(\xi_2)m_N(\xi_3)m_N(\xi_4)}\Big)\\
\times \widehat{\overline{I(|u|^2u)
}}(\xi_1)\widehat{Iu}(\xi_2)\widehat{\overline{Iu}}(\xi_3)\widehat{Iu}(\xi_4)d\xi_2d\xi_3d\xi_4dt\Big|.
\end{split}
\end{equation}
For our purpose, we also adopt a estimate of Coifman-Meyer for a
class of multilinear operators as well as \cite{C04T}. Consider an
infinitely differentiable symbol $m:\R^{nk}\mapsto\C$ so that for
all $\alpha\in \N^{nk}$ and all
$\xi=(\xi_1,\xi_2,\cdots,\xi_k)\in\R^{nk}$, there is a constant
$c(\alpha)$ such that
\begin{align}\label{symbol}
|\partial_{\xi}^\alpha m(\xi)|\leq c(\alpha)(1+|\xi|)^{-|\alpha|}.
\end{align}
Define the multilinear operator $T$ by
\begin{align*}
[T(f_1,\cdots,f_k)](x)=\int_{\R^{nk}}e^{ix\cdot(\xi_1+\cdots+\xi_k)}m(\xi_1,\cdots,\xi_k)\hat
f_1(\xi_1)\cdots\hat f_k(\xi_k)d\xi_1\cdots d\xi_k,
\end{align*}
or
\begin{align*}
\mathcal{F}[T(f_1,\cdots,f_k)](\xi)=\int_{\xi=\xi_1+\cdots+\xi_k}m(\xi_1,\cdots,\xi_k)\hat
f_1(\xi_1)\cdots\hat f_k(\xi_k)d\sigma(\xi).
\end{align*}
\begin{proposition}[\cite{C78M},Page 179.]\label{multiplier}
Suppose $p_j\in (1,\infty), j=1,\cdots k,$ are such that
$\frac1p=\frac1{p_1}+\frac1{p_2}+\cdots+\frac1{p_k}\leq1.$ Assume
$m(\xi_1,\cdots,\xi_k)$ a smooth symbol as in \eqref{symbol}. Then
there is a constant $C=C(p_i,n,k,c(\alpha))$ so that for all Schwarz
class functions $f_1,\cdots,f_k,$
\begin{align}\label{symbol2}
\|[T(f_1,\cdots,f_k)](x)\|_{L^p(\R^n)}\leq
C\|f_1\|_{L^{p_1}(\R^n)}\cdots \|f_k\|_{L^{p_k}(\R^n)}.
\end{align}
\end{proposition}
Now we turn back to our proof of Proposition \ref{ACL}.

{\bf Step 1:} We first estimate $\vartriangle E_1$. To this end, we
decompose
$$u=\sum_{M\geq 1}P_Mu=\sum_{M\geq 1}u_M$$ with the convention that $P_1u:=P_{\leq1}u$.
By utilizing this notation and symmetry, we establish this estimate
\begin{equation}\label{4.17}
\vartriangle E_1\lesssim\sum_{M_1,\cdots, M_4\geq 1\atop{ M_2\geq
M_3\geq M_4}}B(M_1,M_2,M_3,M_4),
\end{equation}
where
\begin{equation}\label{4.18}
\begin{split}
B(M_1,M_2,M_3,M_4):=\Big|\int_0^T\int_{\sum_{j=1}^4\xi_j=0}\Big(1-\frac{m_N(\xi_2+\xi_3+\xi_4)}{m_N(\xi_2)m_N(\xi_3)m_N(\xi_4)}\Big)\\
\times \widehat{\overline{\Delta^2I
u_{M_1}}}(\xi_1)\widehat{Iu_{M_2}}(\xi_2)\widehat{\overline{Iu_{M_3}}}(\xi_3)\widehat{Iu_{M_4}}(\xi_4)d\xi_2d\xi_3d\xi_4dt\Big|.
\end{split}
\end{equation}

{\bf Case $I$}: $M_1>1, M_2\geq M_3\geq M_4>1.$ This case is broken
down into the following several subcases.

{\bf Subcase $I_1$}: $N\gg M_2$. In this case, we have
$$m_N(\xi_2+\xi_3+\xi_4)=m_N(\xi_2)=m_N(\xi_3)=m_N(\xi_4)=1$$
and thus $B(M_1,M_2,M_3,M_4)=0$ and its contribution to the
right-hand side of \eqref{4.17} vanishes.

{\bf Subcase $I_2$}: $M_2\gtrsim N\gg M_3$. Since
$\sum_{j=1}^4\xi_j=0$, we must have $M_1\sim M_2$. The Fundamental
Theorem of Calculus yields that
$$\Big|1-\frac{m_N(\xi_2+\xi_3+\xi_4)}{m_N(\xi_2)m_N(\xi_3)m_N(\xi_4)}\Big|
=\Big|1-\frac{m_N(\xi_2+\xi_3+\xi_4)}{m_N(\xi_2)}\Big|\lesssim
\Big|\frac{\nabla
m_N(\xi_2)\cdot(\xi_3+\xi_4)}{m_N(\xi_2)}\Big|\lesssim
\frac{M_3}{M_2}.$$ By using Proposition \ref{multiplier}, Sobolev
embedding and the Bernstein inequality, and the fact $M_j>1$, we can
see that
\begin{equation*}
\begin{split}
B(M_1,M_2,M_3,M_4)\lesssim \frac{M_3}{M_2}\|\Delta^2I
u_{M_1}\|_{L^4_{t}L^{\frac{2n}{n-2}}_x}\|Iu_{M_2}\|_{L^4_{t}L^{\frac{2n}{n-2}}_x}\|Iu_{M_3}\|_{L^4_{t}L^{n}_x}\|Iu_{M_4}\|_{L^4_{t}L^{n}_x}\\
\lesssim
\frac{M_3M_1^{2}}{M_2^{3}M_3^{\frac{8-n}2}M_4^{\frac{8-n}2}}\|\Delta
I u_{M_1}\|_{L^4_{t}L^{\frac{2n}{n-2}}_x}\|\Delta
Iu_{M_2}\|_{L^4_{t}L^{\frac{2n}{n-2}}_x}\|\Delta
Iu_{M_3}\|_{L^4_{t}L^{\frac{2n}{n-2}}_x}\|\Delta
Iu_{M_4}\|_{L^4_{t}L^{\frac{2n}{n-2}}_x}.
\end{split}
\end{equation*}
and the right hand of above inequality can be controlled by that in
the case for $n=7$
\begin{equation*}
\begin{split}
\frac{1}{M_2}\big(\frac{M_3}{M_4}\big)^{\frac12}Z^4_I(T)\lesssim
N^{-\frac12+}\big(\frac{M_3}{M_2M_4}\big)^{\frac12}M_2^{0-}Z^4_I(T),
\end{split}
\end{equation*}
and in the case for $5\leq n\leq6$
\begin{equation*}
\begin{split}
\frac{1}{M_2}Z^4_I(T)\lesssim N^{-1+}M_2^{0-}Z^4_I(T).
\end{split}
\end{equation*}
The factor $M_2^{0-}$ allows us to sum in $M_1, M_2, M_3, M_4$,
hence we obtain that
\begin{equation}\label{4.19}
\begin{split}
\sum_{M_1,\cdots, M_4\geq 1\atop{ M_2\geq M_3\geq
M_4}}B(M_1,M_2,M_3,M_4)\lesssim
\begin{cases}N^{-\frac12+}Z^4_I(T)\quad\quad n=7;\\ N^{-1+}Z^4_I(T)\quad\quad 5\leq n\leq 6.\end{cases}
\end{split}
\end{equation}

{\bf Subcase $I_3$}: $M_2\gg M_3\gtrsim N$. Since
$\sum_{j=1}^4\xi_j=0$, we must have $M_1\sim M_2$. Observe that
$m_N(\xi_1)\approx m_N(\xi_2)\ll m_N(\xi_3)\lesssim m_N(\xi_4)\leq
1$, it follows that
$$\Big|1-\frac{m_N(\xi_2+\xi_3+\xi_4)}{m_N(\xi_2)m_N(\xi_3)m_N(\xi_4)}\Big|
=\Big|\frac{m_N(\xi_1)-m_N(\xi_2)m_N(\xi_3)m_N(\xi_4)}{m_N(\xi_2)m_N(\xi_3)m_N(\xi_4)}\Big|\lesssim
\frac{m_N(\xi_1)}{m_N(\xi_2)m_N(\xi_3)m_N(\xi_4)}.$$

Applying again the multilinear multiplier theorem, Sobolev embedding
and the Bernstein inequality, and $M_j>1$, we can see that
\begin{equation*}
\begin{split}
&B(M_1,M_2,M_3,M_4)\\&\lesssim
\frac{m_N(\xi_1)}{m_N(\xi_2)m_N(\xi_3)m_N(\xi_4)}\|\Delta^2I
u_{M_1}\|_{L^4_{t}L^{\frac{2n}{n-2}}_x}\|Iu_{M_2}\|_{L^4_{t}L^{\frac{2n}{n-2}}_x}\|Iu_{M_3}\|_{L^4_{t}L^{n}_x}\|Iu_{M_4}\|_{L^4_{t}L^{n}_x}\\
&\lesssim
\frac{m_N(\xi_1)M_3^{\frac{n-8}2}M_4^{\frac{n-8}2}}{m_N(\xi_2)m_N(\xi_3)m_N(\xi_4)}
\|\Delta I u_{M_1}\|_{L^4_{t}L^{\frac{2n}{n-2}}_x}\|\Delta
Iu_{M_2}\|_{L^4_{t}L^{\frac{2n}{n-2}}_x}\|\Delta
Iu_{M_3}\|_{L^4_{t}L^{\frac{2n}{n-2}}_x}\|\Delta
Iu_{M_4}\|_{L^4_{t}L^{\frac{2n}{n-2}}_x},
\end{split}
\end{equation*}
and the right hand of above inequality can be controlled by
\begin{equation*}
\begin{split}
&\frac{M_3^{\frac{n-8}2}M_4^{\frac{n-8}2}}{m_N(\xi_3)m_N(\xi_4)}\|\Delta
I u_{M_1}\|_{L^4_{t}L^{\frac{2n}{n-2}}_x}\|\Delta
Iu_{M_2}\|_{L^4_{t}L^{\frac{2n}{n-2}}_x}Z^2_I(T)\\&\lesssim
\frac{M_3^{\frac{n-8}2}}{m_N(\xi_3)}\|\Delta I
u_{M_1}\|_{L^4_{t}L^{\frac{2n}{n-2}}_x}\|\Delta
Iu_{M_2}\|_{L^4_{t}L^{\frac{2n}{n-2}}_x}Z^2_I(T)\\&\lesssim
N^{\frac{n-8}2+}M_3^{0-}\|\Delta I
u_{M_1}\|_{L^4_{t}L^{\frac{2n}{n-2}}_x}\|\Delta
Iu_{M_2}\|_{L^4_{t}L^{\frac{2n}{n-2}}_x}Z^2_I(T),
\end{split}
\end{equation*}
where we make use of the fact $m_N(\xi)|\xi|^{\frac{8-n}2}$ is
increasing as soon as $s>\frac n2-2$, and the definition of
$m_N(\xi)$. The factor $M_3^{0-}$ allows us to sum in $M_3, M_4$ and
we use the fact $M_1\sim M_2$ and the Cauchy-Schwarz inequality to
get
\begin{equation}\label{4.20}
\begin{split}
&\sum_{M_1,\cdots, M_4\geq 1\atop{ M_2\geq M_3\geq
M_4}}B(M_1,M_2,M_3,M_4)\\&\lesssim
N^{\frac{n-8}2+}\big(\sum_{M_1>1}\|\Delta I
u_{M_1}\|^2_{L^4_{t}L^{\frac{2n}{n-2}}_x}\big)^{\frac12}\big(\sum_{M_2>1}\|\Delta
Iu_{M_2}\|^2_{L^4_{t}L^{\frac{2n}{n-2}}_x}\big)^{\frac12}Z^2_I(T)\\&\lesssim
N^{\frac{n-8}2+}Z^4_I(T).
\end{split}
\end{equation}

{\bf Subcase $I_4$}: $M_2\sim M_3\gtrsim N$. Since
$\sum_{j=1}^4\xi_j=0$, we must have $M_1\lesssim M_2\sim M_3$. A
direct computation yields that
$$\Big|1-\frac{m_N(\xi_2+\xi_3+\xi_4)}{m_N(\xi_2)m_N(\xi_3)m_N(\xi_4)}\Big|
=\Big|\frac{m_N(\xi_1)-m_N(\xi_2)m_N(\xi_3)m_N(\xi_4)}{m_N(\xi_2)m_N(\xi_3)m_N(\xi_4)}\Big|\lesssim
\frac{m_N(\xi_1)}{m_N(\xi_2)m_N(\xi_3)m_N(\xi_4)}.$$

Observe that $m_N(\xi_1)\gtrsim m_N(\xi_2),
m_N(\xi_1)|\xi_1|^2\lesssim m_N(\xi_2)|\xi_2|^2$, thus as similar
argument for case $I_3$ shows that
\begin{equation*}
\begin{split}
&B(M_1,M_2,M_3,M_4)\\&\lesssim
\frac{m_N(\xi_1)}{m_N(\xi_2)m_N(\xi_3)m_N(\xi_4)}\|\Delta^2I
u_{M_1}\|_{L^4_{t}L^{\frac{2n}{n-2}}_x}\|Iu_{M_2}\|_{L^4_{t}L^{\frac{2n}{n-2}}_x}\|Iu_{M_3}\|_{L^4_{t}L^{n}_x}\|Iu_{M_4}\|_{L^4_{t}L^{n}_x}\\
&\lesssim
\frac{m_N(\xi_1)M_1^2M_3^{\frac{n-8}2}M_4^{\frac{n-8}2}}{M_2^2m_N(\xi_2)m_N(\xi_3)m_N(\xi_4)}
\|\Delta I u_{M_1}\|_{L^4_{t}L^{\frac{2n}{n-2}}_x}\|\Delta
Iu_{M_2}\|_{L^4_{t}L^{\frac{2n}{n-2}}_x}\|\Delta
Iu_{M_3}\|_{L^4_{t}L^{\frac{2n}{n-2}}_x}\|\Delta
Iu_{M_4}\|_{L^4_{t}L^{\frac{2n}{n-2}}_x}\\
&\lesssim \frac{M_3^{\frac{n-8}2}}{m_N(\xi_3)} Z^4_I(T)\lesssim
N^{\frac{n-8}2+}M_3^{0-}Z^4_I(T).
\end{split}
\end{equation*}
The factor $M_3^{0-}$ allows us to sum in $M_1,M_2,M_3, M_4$  to
estimate that
\begin{equation}\label{4.21}
\begin{split}
&\sum_{M_1,\cdots, M_4\geq 1\atop{ M_2\geq M_3\geq
M_4}}B(M_1,M_2,M_3,M_4)\lesssim N^{\frac{n-8}2+}Z^4_I(T).
\end{split}
\end{equation}

{\bf Case $II$}: There exits $1\leq j_0\leq 4$ such that
$M_{j_0}=1$. This case is also split into the following several
subcases.

{\bf Subcase $II_1$}: $M_1=1$. In this case, note that $N\gg1$, we
must have $M_2\geq M_3>1=M_4$ or $M_2\geq M_3\geq M_4>1$, otherwise
$$B(M_1,M_2,M_3,M_4)=0.$$
Also, arguing as $I_1$, if $N\gg M_2$ then
$$m_N(\xi_2+\xi_3+\xi_4)=m_N(\xi_2)=m_N(\xi_3)=m_N(\xi_4)=1$$
and thus $B(M_1,M_2,M_3,M_4)=0$ and this contribution to the
right-hand side of \eqref{4.17} vanishes. Therefore, we get
$M_2\gtrsim N$. Furthermore, it follows from $\sum_{j=1}^4\xi_j=0$
that
$$M_2\sim M_3\gtrsim N.$$On the other
hand,  we have
$$\Big|1-\frac{m_N(\xi_2+\xi_3+\xi_4)}{m_N(\xi_2)m_N(\xi_3)m_N(\xi_4)}\Big|
\lesssim
\frac{m_N(\xi_1)}{m_N(\xi_2)m_N(\xi_3)m_N(\xi_4)}=\frac{1}{m_N(\xi_2)m_N(\xi_3)m_N(\xi_4)}.$$
Applying the multilinear multiplier theorem, and Sobolev embedding,
we can see that
\begin{equation*}
\begin{split}
&B(M_1,M_2,M_3,M_4)\\&\lesssim
\frac{1}{m_N(\xi_2)m_N(\xi_3)m_N(\xi_4)}\|\Delta^2I
u_{M_1}\|_{L^4_{t}L^{\frac{2n}{n-2}}_x}\|Iu_{M_2}\|_{L^4_{t}L^{\frac{2n}{n-2}}_x}\|Iu_{M_3}\|_{L^4_{t}L^{n}_x}\|Iu_{M_4}\|_{L^4_{t}L^{n}_x}\\
&\lesssim
\frac{M_1^{2}M_2^{-2}M_3^{\frac{n-8}2}}{m_N(\xi_2)m_N(\xi_3)m_N(\xi_4)}\|\Delta
I u_{M_1}\|_{L^4_{t}L^{\frac{2n}{n-2}}_x}\|\Delta
Iu_{M_2}\|_{L^4_{t}L^{\frac{2n}{n-2}}_x}\|\Delta
Iu_{M_3}\|_{L^4_{t}L^{\frac{2n}{n-2}}_x}\|
Iu_{M_4}\|_{L^4_{t}L^{n}_x}.
\end{split}
\end{equation*}
and the right hand of above inequality can be controlled by, when
$M_2\sim M_3\gtrsim N\gg1=M_4$ as soon as $s>\frac n4-1$,
\begin{equation*}
\begin{split}
\frac{1}{M_2^{2}m_N(\xi_2)m_N(\xi_3)M_3^{\frac{8-n}2}}Z^3_I(T)\|
Iu_{M_4}\|_{L^4_{t}L^{n}_x}&\lesssim
\frac{N^{2(s-2)}}{M_2^{2(1+s-\frac n4)}}Z^3_I(T)\|
Iu_{M_4}\|_{L^4_{t}\dot H^{-\frac{n-5}4}_4}\\&\lesssim N^{-(6-\frac
n2)+}M_2^{0-}Z^3_I(T)\delta,
\end{split}
\end{equation*}
and when $M_2\sim M_3\geq M_4>1$ and $M_2\sim M_3\gtrsim N$ as soon
as $s>\frac n4-1$,
\begin{equation*}
\begin{split}
\frac{M_2^{-2}M_3^{\frac{n-8}2}M_4^{\frac{n-8}2}}{m_N(\xi_2)m_N(\xi_3)m_N(\xi_4)}Z^4_I(T)\lesssim
\frac{M_2^{-2}M_3^{\frac{n-8}2}}{m_N(\xi_2)m_N(\xi_3)}Z^4_I(T)\lesssim
N^{-(6-\frac n2)+}M_2^{0-}Z^4_I(T).
\end{split}
\end{equation*}
The factor $M_2^{0-}$ allows us to sum in $M_1, M_2, M_3, M_4$,
hence we obtain that
\begin{equation}\label{4.22}
\begin{split}
\sum_{M_1,\cdots, M_4\geq 1\atop{ M_2\geq M_3\geq
M_4}}B(M_1,M_2,M_3,M_4)\lesssim N^{-(6-\frac
n2)+}Z^3_I(T)\delta+N^{-(6-\frac n2)+}Z^4_I(T).
\end{split}
\end{equation}

{\bf Subcase $II_2$}: $M_1>1$. This subcase is split into the
following several subcases.

{\bf Sub-subcase $II_2^a$}: $M_1>1, M_2=M_3=M_4=1$. In this case,
since $\sum\limits_{j=1}^4\xi_j=0$, we must have $M_1\sim 1\ll N$.
Then again we have
$$m_N(\xi_2+\xi_3+\xi_4)=m_N(\xi_2)=m_N(\xi_3)=m_N(\xi_4)=1$$
and thus $B(M_1,M_2,M_3,M_4)=0$.

{\bf Sub-subcase $II_2^b$}: $M_1>1, M_2>1=M_3=M_4$. In this case, we
must have $M_1\sim M_2$ since $\sum_{j=1}^4\xi_j=0$. We may assume
that $M_1\sim M_2\gtrsim N$, since otherwise $B(M_1,M_2,M_3,M_4)=0$.
Now our purpose is to estimate $B(M_1,M_2,M_3,M_4)$ under the
circumstance that
$$M_1\sim M_2\gtrsim N\gg1=M_3=M_4.$$
In addition,  the Fundamental Theorem of Calculus yields
$$\Big|1-\frac{m_N(\xi_2+\xi_3+\xi_4)}{m_N(\xi_2)m_N(\xi_3)m_N(\xi_4)}\Big|
=\Big|1-\frac{m_N(\xi_2+\xi_3+\xi_4)}{m_N(\xi_2)}\Big|\lesssim
\Big|\frac{\nabla
m_N(\xi_2)\cdot(\xi_3+\xi_4)}{m_N(\xi_2)}\Big|\lesssim
\frac{1}{M_2}.$$ Applying Proposition \ref{multiplier}, Sobolev
embedding, and the Bernstein inequality, we can see that
\begin{equation*}
\begin{split}
B(M_1,M_2,M_3,M_4)\lesssim \frac{1}{M_2}\|\Delta^2I
u_{M_1}\|_{L^4_{t}L^{\frac{2n}{n-2}}_x}\|Iu_{M_2}\|_{L^4_{t}L^{\frac{2n}{n-2}}_x}\|Iu_{M_3}\|_{L^4_{t}L^{n}_x}\|Iu_{M_4}\|_{L^4_{t}L^{n}_x}\\
\lesssim \frac{M_1^{2}}{M_2^{3}}\|\Delta I
u_{M_1}\|_{L^4_{t}L^{\frac{2n}{n-2}}_x}\|\Delta
Iu_{M_2}\|_{L^4_{t}L^{\frac{2n}{n-2}}_x}\|Iu_{M_3}\|_{L^4_{t}L^{n}_x}\|Iu_{M_4}\|_{L^4_{t}L^{n}_x}.
\end{split}
\end{equation*}
and the right hand of above inequality can be controlled by
\begin{equation*}
\begin{split}
\frac{1}{M_2}Z^2_I(T)\|Iu_{M_3}\|_{L^4_{t}\dot
H^{-\frac{n-5}{4}}_4}\|Iu_{M_4}\|_{L^4_{t}\dot
H^{-\frac{n-5}{4}}_4}\lesssim N^{-1+}M_2^{0-}Z^2_I(T)\delta^2,
\end{split}
\end{equation*}
The factor $M_2^{0-}$ allows us to sum in $M_1, M_2, M_3, M_4$,
hence we obtain that
\begin{equation}\label{4.23}
\begin{split}
\sum_{M_1,\cdots, M_4\geq 1\atop{ M_2\geq M_3\geq
M_4}}B(M_1,M_2,M_3,M_4)\lesssim N^{-1+}Z^2_I(T)\delta^2.
\end{split}
\end{equation}

{\bf Sub-subcase $II_2^c$}: $M_1>1, M_2\geq M_3>1=M_4$. We may
assume $M_2\gtrsim N$, otherwise the contribution of this case is
null. Arguing similar as for Case $I$, we also break this case into
several cases.

{\bf $\clubsuit$} If $M_2\gtrsim N\gg M_3>1=M_4$, we must have
$M_1\sim M_2\gtrsim N\gg M_3>1=M_4$ since $\sum_{j=1}^4\xi_j=0$.
Hence, we obtain that
$$\Big|1-\frac{m_N(\xi_2+\xi_3+\xi_4)}{m_N(\xi_2)m_N(\xi_3)m_N(\xi_4)}\Big|
=\Big|1-\frac{m_N(\xi_2+\xi_3+\xi_4)}{m_N(\xi_2)}\Big|\lesssim
\Big|\frac{\nabla
m_N(\xi_2)\cdot(\xi_3+\xi_4)}{m_N(\xi_2)}\Big|\lesssim
\frac{M_3}{M_2}.$$ Applying the multilinear multiplier theorem,
Sobolev embedding, and the Bernstein inequality and keeping in mind
$M_1, M_2, M_3>1, M_4=1$, we can see that
\begin{equation*}
\begin{split}
B(M_1,M_2,M_3,M_4)\lesssim \frac{M_3}{M_2}\|\Delta^2I
u_{M_1}\|_{L^4_{t}L^{\frac{2n}{n-2}}_x}\|Iu_{M_2}\|_{L^4_{t}L^{\frac{2n}{n-2}}_x}\|Iu_{M_3}\|_{L^4_{t}L^{n}_x}\|Iu_{M_4}\|_{L^4_{t}L^{n}_x}\\
\lesssim \frac{M_3M_1^{2}}{M_2^{3}M_3^{\frac{8-n}2}}\|\Delta I
u_{M_1}\|_{L^4_{t}L^{\frac{2n}{n-2}}_x}\|\Delta
Iu_{M_2}\|_{L^4_{t}L^{\frac{2n}{n-2}}_x}\|\Delta
Iu_{M_3}\|_{L^4_{t}L^{\frac{2n}{n-2}}_x}\| Iu_{M_4}\|_{L^4_{t}\dot
H^{-\frac{n-5}{4}}_4}.
\end{split}
\end{equation*}
and the right hand of above inequality can be controlled by that in
the case that $n=7$
\begin{equation*}
\begin{split}
\frac{{M_3}^{\frac12}}{M_2}Z^3_I(T)\delta\lesssim
N^{-\frac12+}M_2^{0-}Z^3_I(T)\delta,
\end{split}
\end{equation*}
and in the case that $5\leq n\leq6$
\begin{equation*}
\begin{split}
\frac{1}{M_2}Z^3_I(T)\delta\lesssim N^{-1+}M_2^{0-}Z^3_I(T)\delta.
\end{split}
\end{equation*}
The factor $M_2^{0-}$ allows us to sum in $M_1, M_2, M_3, M_4$,
hence we obtain that
\begin{equation}\label{4.24}
\begin{split}
\sum_{M_1,\cdots, M_4\geq 1\atop{ M_2\geq M_3\geq
M_4}}B(M_1,M_2,M_3,M_4)\lesssim
\begin{cases}N^{-\frac12+}Z^3_I(T)\delta\quad\quad n=7;\\ N^{-1+}Z^3_I(T)\delta\quad\quad 5\leq n\leq 6.\end{cases}
\end{split}
\end{equation}

{\bf $\clubsuit$} $M_2\gg M_3\gtrsim N$. Since
$\sum_{j=1}^4\xi_j=0$, we must have $M_1\sim M_2$. Observe that
$m_N(\xi_1)\approx m_N(\xi_2)\ll m_N(\xi_3)\lesssim m_N(\xi_4)=1$,
it follows that
$$\Big|1-\frac{m_N(\xi_2+\xi_3+\xi_4)}{m_N(\xi_2)m_N(\xi_3)m_N(\xi_4)}\Big|
=\Big|\frac{m_N(\xi_1)-m_N(\xi_2)m_N(\xi_3)m_N(\xi_4)}{m_N(\xi_2)m_N(\xi_3)}\Big|\lesssim
\frac{m_N(\xi_1)}{m_N(\xi_2)m_N(\xi_3)}.$$ Applying again the
multilinear multiplier theorem, Sobolev embedding and the Bernstein
inequality and $M_j>1$ for $j=1,2,3$ and $M_4=1$, we can see that
\begin{equation*}
\begin{split}
&B(M_1,M_2,M_3,M_4)\\&\lesssim
\frac{m_N(\xi_1)}{m_N(\xi_2)m_N(\xi_3)m_N(\xi_4)}\|\Delta^2I
u_{M_1}\|_{L^4_{t}L^{\frac{2n}{n-2}}_x}\|Iu_{M_2}\|_{L^4_{t}L^{\frac{2n}{n-2}}_x}\|Iu_{M_3}\|_{L^4_{t}L^{n}_x}\|Iu_{M_4}\|_{L^4_{t}L^{n}_x}\\
&\lesssim \frac{m_N(\xi_1)M_3^{\frac{n-8}2}}{m_N(\xi_2)m_N(\xi_3)}
\|\Delta I u_{M_1}\|_{L^4_{t}L^{\frac{2n}{n-2}}_x}\|\Delta
Iu_{M_2}\|_{L^4_{t}L^{\frac{2n}{n-2}}_x}\|\Delta
Iu_{M_3}\|_{L^4_{t}L^{\frac{2n}{n-2}}_x}\| Iu_{M_4}\|_{L^4_{t}\dot
H^{-\frac{n-5}{4}}_4}.
\end{split}
\end{equation*}
and the right hand of above inequality can be controlled by
\begin{equation*}
\begin{split}
&\frac{M_3^{\frac{n-8}2}}{m_N(\xi_3)}\|\Delta I
u_{M_1}\|_{L^4_{t}L^{\frac{2n}{n-2}}_x}\|\Delta
Iu_{M_2}\|_{L^4_{t}L^{\frac{2n}{n-2}}_x}Z_I(T)\delta\\&\lesssim
\frac{M_3^{\frac{n-8}2}}{m_N(\xi_3)}\|\Delta I
u_{M_1}\|_{L^4_{t}L^{\frac{2n}{n-2}}_x}\|\Delta
Iu_{M_2}\|_{L^4_{t}L^{\frac{2n}{n-2}}_x}Z_I(T)\delta\\&\lesssim
N^{\frac{n-8}2+}M_3^{0-}\|\Delta I
u_{M_1}\|_{L^4_{t}L^{\frac{2n}{n-2}}_x}\|\Delta
Iu_{M_2}\|_{L^4_{t}L^{\frac{2n}{n-2}}_x}Z_I(T)\delta.
\end{split}
\end{equation*}
The factor $M_3^{0-}$ allows us to sum in $M_3, M_4$ and the fact
$M_1\sim M_2$ and the Cauchy-Schwarz inequality permit us to
estimate that
\begin{equation}\label{4.25}
\begin{split}
&\sum_{M_1,\cdots, M_4\geq 1\atop{ M_2\geq M_3\geq
M_4}}B(M_1,M_2,M_3,M_4)\\&\lesssim
N^{\frac{n-8}2+}\big(\sum_{M_1>1}\|\Delta I
u_{M_1}\|^2_{L^4_{t}L^{\frac{2n}{n-2}}_x}\big)^{\frac12}\big(\sum_{M_2>1}\|\Delta
Iu_{M_2}\|^2_{L^4_{t}L^{\frac{2n}{n-2}}_x}\big)^{\frac12}Z_I(T)\delta\\&\lesssim
N^{\frac{n-8}2+}Z^3_I(T)\delta.
\end{split}
\end{equation}

{\bf $\clubsuit$} $M_2\sim M_3\gtrsim N$. Since
$\sum_{j=1}^4\xi_j=0$, we must have $M_1\lesssim M_2\sim M_3$. A
direct computation yields that
$$\Big|1-\frac{m_N(\xi_2+\xi_3+\xi_4)}{m_N(\xi_2)m_N(\xi_3)m_N(\xi_4)}\Big|
=\Big|\frac{m_N(\xi_1)-m_N(\xi_2)m_N(\xi_3)m_N(\xi_4)}{m_N(\xi_2)m_N(\xi_3)m_N(\xi_4)}\Big|\lesssim
\frac{m_N(\xi_1)}{m_N(\xi_2)m_N(\xi_3)}.$$

Observe that $m_N(\xi_1)\gtrsim m_N(\xi_2),
m_N(\xi_1)|\xi_1|^2\lesssim m_N(\xi_2)|\xi_2|^2$, thus similar
argument as above leads to that
\begin{equation*}
\begin{split}
&B(M_1,M_2,M_3,M_4)\\&\lesssim
\frac{m_N(\xi_1)}{m_N(\xi_2)m_N(\xi_3)m_N(\xi_4)}\|\Delta^2I
u_{M_1}\|_{L^4_{t}L^{\frac{2n}{n-2}}_x}\|Iu_{M_2}\|_{L^4_{t}L^{\frac{2n}{n-2}}_x}\|Iu_{M_3}\|_{L^4_{t}L^{n}_x}\|Iu_{M_4}\|_{L^4_{t}L^{n}_x}\\
&\lesssim
\frac{m_N(\xi_1)M_1^2M_3^{\frac{n-8}2}}{M_2^2m_N(\xi_2)m_N(\xi_3)}
\|\Delta I u_{M_1}\|_{L^4_{t}L^{\frac{2n}{n-2}}_x}\|\Delta
Iu_{M_2}\|_{L^4_{t}L^{\frac{2n}{n-2}}_x}\|\Delta
Iu_{M_3}\|_{L^4_{t}L^{\frac{2n}{n-2}}_x}\| Iu_{M_4}\|_{L^4_{t}\dot
H^{-\frac{n-5}{4}}_4}\\
&\lesssim \frac{M_3^{\frac{n-8}2}}{m_N(\xi_3)}
Z^3_I(T)\delta\lesssim N^{\frac{n-8}2+}M_3^{0-}Z^3_I(T)\delta.
\end{split}
\end{equation*}
The factor $M_3^{0-}$ allows us to sum in $M_1,M_2,M_3, M_4$  to
estimate that
\begin{equation}\label{4.26}
\begin{split}
&\sum_{M_1,\cdots, M_4\geq 1\atop{ M_2\geq M_3\geq
M_4}}B(M_1,M_2,M_3,M_4)\lesssim N^{\frac{n-8}2+}Z^3_I(T)\delta.
\end{split}
\end{equation}
Putting all of cases together, it follows from
\eqref{4.19}-\eqref{4.26} that
\begin{equation}\label{4.27}
\begin{split}
\vartriangle E_1\lesssim \max\{N^{-1+},
N^{\frac{n-8}2+}\}\Big(Z^4_I(T)+Z^3_I(T)\delta+Z^2_I(T)\delta^2\Big).
\end{split}
\end{equation}

{\bf Step 2:} We secondly estimate $\vartriangle E_2$. To this end,
we again decompose
$$u=\sum_{M\geq 1}P_Mu=\sum_{M\geq 1}u_M$$ with the convention that $P_1u:=P_{\leq1}u$.
By utilizing this notation and symmetry, we establish this estimate
\begin{equation}\label{4.28}
\vartriangle E_2\lesssim\sum_{M_1,\cdots, M_4\geq 1\atop{ M_2\geq
M_3\geq M_4}}C(M_1,M_2,M_3,M_4)
\end{equation}
where
\begin{equation}\label{4.29}
\begin{split}
C(M_1,M_2,M_3,M_4):=\Big|\int_0^T\int_{\sum_{j=1}^4\xi_j=0}\Big(1-\frac{m_N(\xi_2+\xi_3+\xi_4)}{m_N(\xi_2)m_N(\xi_3)m_N(\xi_4)}\Big)\\
\times \widehat{\overline{P_{M_1}I
(|u|^2u)}}(\xi_1)\widehat{Iu_{M_2}}(\xi_2)\widehat{\overline{Iu_{M_3}}}(\xi_3)\widehat{Iu_{M_4}}(\xi_4)d\xi_2d\xi_3d\xi_4dt\Big|.
\end{split}
\end{equation}
In order to estimate $C(M_1,M_2,M_3,M_4)$, we make the observation
that in estimating $B(M_1,M_2,M_3,M_4)$ for the term involving the
$M_1$ frequency we only used the bound
\begin{equation}\label{4.30}
\|\Delta^2I u_{M_1}\|_{L^4_{t}L^{\frac{2n}{n-2}}_x}\leq M_1^2
\|\Delta I u_{M_1}\|_{L^4_{t}L^{\frac{2n}{n-2}}_x}\lesssim M_1^2
Z_I(t).
\end{equation}
Thus to estimate $\vartriangle E_2$, it suffices to show that
\begin{equation}\label{4.31}
\|P_{M_1}I (|u|^2u)\|_{L^4_{t}L^{\frac{2n}{n-2}}_x}\lesssim
M_1^2Z^3_I(t)
\end{equation}
and then arguing as for estimating $\vartriangle E_1$, we substitute
\eqref{4.31} for \eqref{4.30} to obtain that
\begin{equation}\label{4.32}
\begin{split}
\vartriangle E_2\lesssim \max\{N^{-1+},
N^{\frac{n-8}2+}\}\Big(Z^6_I(T)+Z^5_I(T)\delta+Z^4_I(T)\delta^2\Big).
\end{split}
\end{equation}
Therefore, we are left to prove \eqref{4.31}. The boundedness of the
Littlewood-Paley operator and the Sobolev embedding yield that
\begin{equation*}
\begin{split}
M_1^{-2}\|P_{M_1}I
(|u|^2u)\|_{L^4_{t}L^{\frac{2n}{n-2}}_x}\lesssim&\|P_{M_1}I
(|u|^2u)\|_{L^4_{t}L^{\frac{2n}{n+2}}_x}\leq
\|u\|^3_{L^{12}_{t}L^{\frac{6n}{n+2}}_x}.
\end{split}
\end{equation*}
We decompose $u$ into low frequency and high frequency like that
$u:=u_{\leq N}+u_{> N}$. We first estimate the low frequency part by
interpolation
\begin{equation}\label{4.33}
\| u_{\leq N}\|^3_{L^{12}_{t}L^{\frac{6n}{n+2}}_x}\leq \| Iu_{\leq
N}\|^3_{L^{12}_{t}L^{\frac{6n}{n+2}}_x}\leq \| Iu_{\leq
N}\|^{3(1-\theta)}_{L^{12}_{t}L^{\frac{6n}{3n-2}}_x}\| Iu_{\leq
N}\|^{3\theta}_{L^{12}_{t}L^{\frac{6n}{3n-14}}_x}\leq Z^3_I(t),
\end{equation}
with $\theta=\frac{n-2}6$. For the high frequency, we have that
\begin{equation}\label{4.34}
\| u_{>N}\|^3_{L^{12}_{t}L^{\frac{6n}{n+2}}_x}\leq \|
|\nabla|^{\frac{n-2}3}u_{>
N}\|^3_{L^{12}_{t}L^{\frac{6n}{3n-2}}_x}=\|
|\nabla|^{\frac{n-2}3}(\Delta I)^{-1}\Delta Iu_{>
N}\|^3_{L^{12}_{t}L^{\frac{6n}{3n-2}}_x}.
\end{equation}
We also can rewrite the right hand as follows
\begin{equation*}
N^{\frac{n-8}3}\|\mathcal{F}^{-1} (\sigma(\xi))\Delta Iu_{>
N}\|^3_{L^{12}_{t}L^{\frac{6n}{3n-2}}_x}
\end{equation*}
with $\sigma(\xi)=(N|\xi|^{-1})^{s-\frac{n-2}3}$ and it follows from
$s\geq\frac{n-2}3$ and H\"ormander's multiplier theorem that
\begin{equation*}
N^{\frac{n-8}3}\|\mathcal{F}^{-1} (\sigma(\xi))\Delta Iu_{>
N}\|^3_{L^{12}_{t}L^{\frac{6n}{3n-2}}_x}\leq N^{\frac{n-8}3}\|\Delta
Iu_{> N}\|^3_{L^{12}_{t}L^{\frac{6n}{3n-2}}_x}\leq Z^3_I(t),
\end{equation*}
since $N\gg1$. This together with \eqref{4.34} gives that
\begin{equation}\label{4.35}
\| u_{>N}\|^3_{L^{12}_{t}L^{\frac{6n}{n+2}}_x}\lesssim Z^3_I(t).
\end{equation}
Finally, \eqref{4.31} follows from \eqref{4.33} and \eqref{4.35} and
this completes the proof of the almost conservation law Proposition
\ref{ACL}.
\end{proof}
\section{Proof of Main Theorem}
We combine the interaction Morawetz estimate and almost conservation
law with a scaling argument to prove the following statement giving
uniform bounds in terms of the rough norm of the initial data.
\begin{proposition}
Suppose $u(t,x)$ is a global in time solution to \eqref{Eq} from
data $u_0\in C_0^\infty(\R^n)$. Then so long as $s> s_0$ with $s_0$
in Theorem \ref{the1}, we have
\begin{equation}\label{5.1}
\|u\|_{M(\R)}\leq C(\|u_0\|_{H^s(\R^n)}),
\end{equation}
\begin{equation}\label{5.2}
\sup_{0\leq t<\infty}\|u\|_{H^s(\R^n)}\leq C(\|u_0\|_{H^s(\R^n)}).
\end{equation}
\end{proposition}
\begin{remark}
The global well-posedness part of Theorem \ref{the1} follows from
\eqref{5.2}, Proposition \ref{Local}  and the standard density
argument.
\end{remark}
\begin{proof}
If $u$ is a solution to \eqref{Eq}, then so is
\begin{equation}\label{5.3}
u^\lambda(t,x)=\lambda^{-2}u(\lambda^{-4}t,\lambda^{-1}x).
\end{equation}
We choose $\lambda$ so that $E(Iu_0^\lambda)=\frac12\|\Delta
Iu_0^\lambda\|_{L^2(\R^n)}^2+\frac14\|Iu_0^\lambda\|_{L^4(\R^n)}^4\lesssim
1$ to remove the uniform bound condition $E(Iu_0)$ in Proposition
\ref{ACL}. As in \eqref{4.2}, we show
\begin{equation}\label{5.6}
\|\Delta Iu_0^\lambda\|_{L^2(\R^n)}\lesssim N^{2-s}\lambda^{\frac
n2-2-s}\|u_0\|_{H^s(\R^n)},
\end{equation}
then the right choice of $\lambda$ is
\begin{equation}\label{5.4}
\lambda\approx N^{\frac{2-s}{s-(\frac n2-2)}}.
\end{equation}
We estimate $\|Iu_0^\lambda\|_{L^4(\R^n)}$, the second term in
$E(Iu_0^\lambda)$, by separating the domains in the frequency space.
Set
\begin{equation*}
\widehat{u_0^{\lambda}}(\xi)=\big(\chi_0(\xi)+\chi_1(\xi)+\chi_2(\xi)\big)\widehat{u_0^{\lambda}}(\xi),
\end{equation*}
for nonnegative smooth functions $\chi_j(\xi)$ such that
$\sum_{j=0}^2\chi_j(\xi)=1$ and $\chi_j$ is supported in
$\{\xi:|\xi|\leq\frac2\lambda\}$, $\{\xi:\frac1\lambda\leq |\xi|\leq
N\}$ and $\{\xi:|\xi|\geq \frac N2\}$, respectively. Then
\begin{equation*}
\widehat{Iu_0^{\lambda}}(\xi)=
\chi_0(\xi)\widehat{u_0^{\lambda}}(\xi)+\chi_1(\xi)\widehat{u_0^{\lambda}}(\xi)+\chi_2(\xi)m_N(\xi)\widehat{u_0^{\lambda}}(\xi).
\end{equation*}
A straightforward argument using Sobolev embedding together with the
relation \eqref{5.4} will give
\begin{equation*}
\|\mathcal{F}^{-1}\big(\chi_0(\xi)\widehat{u_0^{\lambda}}(\xi)\big)\|_{L^4(\R^n)}\lesssim
\lambda^{\frac n4-2}\|u_0\|_{L^2(\R^n)}.
\end{equation*}
\begin{equation*}
\begin{split}
\|\mathcal{F}^{-1}\big(m_N(\xi)\chi_2(\xi)\widehat{u_0^{\lambda}}(\xi)\big)\|_{L^4(\R^n)}&\lesssim
\|\big(\frac{N}{|\xi|}\big)^{2-s}|\xi|^{\frac n
4-s}|\xi|^{s}\chi_2(\xi)\widehat{u_0^{\lambda}}(\xi)\|_{L^2(\R^n)}\\&\lesssim
N^{\frac n4-s}\lambda^{-(s+2-\frac n2)}\|u_0\|_{H^s(\R^n)}.
\end{split}
\end{equation*}
For the medium frequency, we similarly have that
\begin{equation*}
\begin{split}
\|\mathcal{F}^{-1}\big(\chi_1(\xi)\widehat{u_0^{\lambda}}(\xi)\big)\|_{L^4(\R^n)}&\lesssim
\||\xi|^{\frac n
4-s}|\xi|^{s}\chi_1(\xi)\widehat{u_0^{\lambda}}(\xi)\|_{L^2(\R^n)}\\&\leq\begin{cases}N^{\frac
n4-s}\lambda^{-(2+s-\frac n2)}\|u_0\|_{\dot H^s(\R^n)}\quad s\le
\frac n4
\\ \lambda^{\frac n4-2}\|u_0\|_{\dot H^s(\R^n)}\quad s>\frac n4\end{cases} .
\end{split}
\end{equation*}
Summing up the three parts, we obtain that by \eqref{5.4}
\begin{equation}\label{5.5}
\begin{split}
\|Iu_0^\lambda\|_{L^4(\R^n)}&\lesssim \big(\lambda^{\frac
n4-2}+\lambda^{\frac{(s+2-\frac n2)(\frac
n4-2)}{2-s}}\big)\|u_0\|_{H^s(\R^n)}.
\end{split}
\end{equation}
Thus, taking $\lambda$ sufficiently large depending on
$\|u_0\|_{H^s}$ and $N$ (which will be chosen later and will depend
only on $\|u_0\|_{H^s}$), it follows from \eqref{5.6} and
\eqref{5.5} that
\begin{equation}\label{5.7}
\begin{split}
E(Iu_0^\lambda)\lesssim 1.
\end{split}
\end{equation}
We now show that there exists an absolute constant $C_1$ such that
\begin{equation}\label{5.8}
\|u^\lambda\|_{M(\R)}\leq C_1\lambda^{\frac74(\frac n 4-1)}.
\end{equation}
Undoing the scaling, this yields \eqref{5.1}. We prove \eqref{5.8}
via a bootstrap argument. By time reversal symmetry, it suffices to
argue for positive times only. Define
\begin{equation*}
\Omega_1:=\{t\in[0,\infty):\|u^\lambda\|_{M([0,t])}\leq
C_1\lambda^{\frac74(\frac n 4-1)}\}.
\end{equation*}
We want to show $\Omega_1=[0,\infty).$ Let
\begin{equation*}
\Omega_2:=\{t\in[0,\infty):\|u^\lambda\|_{M([0,t])}\leq
2C_1\lambda^{\frac74(\frac n 4-1)}\}.
\end{equation*}
In order to run the bootstrap argument successfully, we need to
verify four things:

1) $\Omega_1\neq \emptyset$. This is obvious as $0 \in \Omega_1$.

2) $\Omega_1$ is closed. This follows from Fatou's Lemma.

3) $\Omega_2\subset \Omega_1$.

4) If $T\in \Omega_1$, then there exists $\varepsilon>0$ such that
$[T, T+\varepsilon)\subset \Omega_1$. This is a consequence of the
local well-posedness theory and 3). We skip the details.

Thus, we need to prove 3). Fix $T\in \Omega_2$; we will show that
$T\in \Omega_1$. By the interaction Morawetz estimate \eqref{IMorE}
and the mass conservation, we can see that
\begin{equation}\label{5.9}
\|u^\lambda\|_{M([0,T])}\leq
\|u_0^\lambda\|_{L^2}^{\frac12}\|u^\lambda\|_{L^\infty([0,T];\dot
H^{\frac12}(\R^n))}^{\frac12}\lesssim_{\|u_0\|_{L^2}} \lambda^{\frac
n4-1}\|u^\lambda\|_{L^\infty([0,T];\dot
H^{\frac12}(\R^n))}^{\frac12}.
\end{equation}
To control the second factor $\|u^\lambda\|_{L^\infty([0,T];\dot
H^{\frac12}(\R^n))}$, we decompose
\begin{equation*}
u^\lambda(t)=P_{\leq N}u^\lambda(t)+P_{>N}u^\lambda(t).
\end{equation*}
In order to estimate the low frequencies, we interpolate between the
$L^2_x$-norm and $\dot H^2_x$-norm and use the fact that the
operator $I$ is the identity on frequencies $|\xi|\leq N$:
\begin{equation}\label{5.10}
\|P_{\leq N}u^\lambda(t)\|_{\dot H^{\frac12}_x}\lesssim
\|u^\lambda(t)\|_{L^2}^{\frac34}\|Iu^\lambda(t)\|_{\dot
H^{2}}^{\frac14}\lesssim_{\|u_0\|_{L^2}} \lambda^{\frac
{3n}8-\frac32}\|Iu^\lambda(t)\|_{\dot H^{2}}^{\frac14}.
\end{equation}
To dominate the high frequencies, we interpolate between the
$L^2_x$-norm and $\dot H^s_x$-norm and use the definition of
operator $I$ to get:
\begin{equation}\label{5.11}
\begin{split}
\|P_{> N}u^\lambda(t)\|_{\dot H^{\frac12}_x}\lesssim
\|u^\lambda(t)\|_{L^2}^{1-\frac1{2s}}\|P_{>N}u^\lambda(t)\|_{\dot
H^{s}}^{\frac1{2s}}&\lesssim_{\|u_0\|_{L^2}} \lambda^{(\frac
{n}2-2)(1-\frac1{2s})}N^{\frac{s-2}{2s}}\|Iu^\lambda(t)\|_{\dot
H^{2}}^{\frac1{2s}}\\ &\lesssim_{\|u_0\|_{L^2}} \lambda^{\frac
{n}2-\frac52}\|Iu^\lambda(t)\|_{\dot H^{2}}^{\frac1{2s}}.
\end{split}
\end{equation}
Collecting \eqref{5.10} through \eqref{5.11}, we obtain that
\begin{equation}\label{5.12}
\begin{split}
\|u^\lambda(t)\|_{M([0,T])}&\lesssim \lambda^{\frac
n4-1}\sup_{t\in[0,T]}\big(\lambda^{\frac
{3n}{16}-\frac34}\|Iu^\lambda(t)\|_{\dot
H^{2}}^{\frac18}+\lambda^{\frac {n}4-\frac54}\|Iu^\lambda(t)\|_{\dot
H^{2}}^{\frac1{4s}}\big)\\&\lesssim \lambda^{\frac74 (\frac
n4-1)}\sup_{t\in[0,T]}\big(\|Iu^\lambda(t)\|_{\dot
H^{2}}^{\frac18}+\|Iu^\lambda(t)\|_{\dot H^{2}}^{\frac1{4s}}\big),
\end{split}
\end{equation}
where we make use of the facts that $\lambda\gg 1$ and $n< 8$ in the
last inequality. Thus, taking $C_1$ sufficiently large depending on
$\|u_0\|_{L^2_x}$, we obtain $T\in \Omega_1$, provided that
\begin{equation}\label{5.13}
\begin{split}
\sup_{t\in[0,T]}\|Iu^\lambda(t)\|_{\dot H^{2}}\leq 1 .
\end{split}
\end{equation}
We now prove that \eqref{5.13} when $T\in\Omega_2$. In practice, let
$\delta>0$ be sufficiently small constant as in Proposition
\ref{ACL}, and we divide $[0,T]$ into
\begin{equation}\label{5.14}
\begin{split}
L\sim \Big(\frac{\lambda^{\frac74(\frac n4-1)}}\delta\Big)^4
\end{split}
\end{equation}
sub-intervals $I_j=[t_j,t_{j+1}]$ such that
\begin{equation*}
\begin{split}
\|u^\lambda\|_{M(I_j)}\leq \delta.
\end{split}
\end{equation*}
Applying Proposition \ref{ACL} on each of the sub-intervals $I_j$,
we get that
\begin{equation*}
\begin{split}
\sup_{t\in[0,T]}E(Iu^\lambda(t))\leq
E(Iu_0^\lambda)+LN^{\max\{-1,\frac{n-8}2\}+}.
\end{split}
\end{equation*}
To maintain small energy during the iteration, we need
\begin{equation*}
\begin{split}
LN^{\max\{-1,\frac{n-8}2\}+}\sim \lambda^{7(\frac
n4-1)}N^{\max\{-1,\frac{n-8}2\}+}\delta^{-4}\ll1,
\end{split}
\end{equation*}
which combined with \eqref{5.4} leads to
\begin{equation*}
\begin{split}
\bigg(N^{\frac{2-s}{s+2-\frac n2}}\bigg)^{7(\frac
n4-1)}N^{\max\{-1,\frac{n-8}2\}+}\ll1.
\end{split}
\end{equation*}
This may be ensured by taking $N=N(\|u_0\|_{H^s})$ large enough
provided that
\begin{equation}\label{5.15}
s>~\begin{cases} \frac{16(n-4)}{7n-24}\qquad5\leq n\leq 6,\\
\frac{45}{23}\qquad\qquad n=7.
\end{cases}
\end{equation}
This completes the bootstrap argument and hence \eqref{5.8} and
moreover \eqref{5.1} follows. To estimate $\|u(T)\|_{H^s_x}$, we
write that by the conservation of mass and the scaling
\begin{equation*}
\begin{split}
\|u(T)\|_{H^s_x}\lesssim \|u_0\|_{L^2_x}+\|u(T)\|_{\dot H^s_x}
\lesssim \|u_0\|_{L^2_x}+\lambda^{s+2-\frac
n2}\|u^\lambda(\lambda^4T)\|_{\dot H^s_x}.
\end{split}
\end{equation*}
Utilizing \eqref{4.3}, the right hand can be controlled by
\begin{equation*}
\begin{split}
\|u_0\|_{L^2_x}+\lambda^{s+2-\frac n2}\|Iu^\lambda(\lambda^4T)\|_{
H^2_x}\lesssim \|u_0\|_{L^2_x}+\lambda^{s+2-\frac
n2}\big(\|u^\lambda(\lambda^4T)\|_{
L^2_x}+\|Iu^\lambda(\lambda^4T)\|_{ \dot H^2_x}\big)
\end{split}
\end{equation*}
Therefore, it follows from \eqref{5.13} that for all $T\in \R$
\begin{equation*}
\begin{split}
\|u(T)\|_{H^s_x} \lesssim \|u_0\|_{L^2_x}+\lambda^{s+2-\frac
n2}\big(\lambda^{\frac n2-2}\|u_0\|_{ L^2_x}+1\big) \lesssim
C(\|u_0\|_{H^s_x}).
\end{split}
\end{equation*}
Hence, we have
\begin{equation}\label{5.16}
\|u(t)\|_{L^\infty(\R;H^s_x)}\lesssim C(\|u_0\|_{H^s_x}).
\end{equation}

{\bf Scattering}\quad We prove that scattering holds in $H^s_x$ for
$s>s_0$. We first show that the global Morawetz estimate \eqref{5.1}
can be upgraded to the global Strichartz estimate
\begin{equation}\label{5.17}
\|u\|_{S^s(I)}:=\sup_{(q_0,r_0) \in
\Lambda_0}\||\nabla|^{s+\frac2{q_0}}u\|_{L^{q_0}_tL^{r_0}_x(I\times\R^n)}+\sup_{(q_1,r_1)
\in \Lambda_1}\|u\|_{L^{q_1}_tL^{r_1}_x(I\times\R^n)}.
\end{equation}
The second step is to use this estimate to prove asymptotic
completeness. The construction of the wave operator is a standard
step, which we omit it here.

Let $u$ be a global solution to \eqref{Eq} with initial data in
$H^s(\R^n)$ for $s>s_0$. From the global Morawetz estimate
\eqref{5.1}, we have
\begin{equation*}
\|u\|_{M(\R)}\leq C(\|u_0\|_{H^s_x}).
\end{equation*}
Let $\delta>0$ be a small constant to be chosen momentarily and
split $\R$ into $L=L(\|u_0\|_{H^s_x})$ subintervals
$I_j=[t_j,t_{j+1}]$ such that
\begin{equation}\label{5.18}
\|u\|_{M(I_j)}\leq \delta.
\end{equation}
By the Strichartz estimate in Proposition \ref{Strichartz}, we have
\begin{equation}\label{5.19}
\|u\|_{S^s(I_j)}\leq\|\langle\nabla\rangle^su(t_j)\|_{L^2_x}+
\||\nabla|^{s-1}(|u|^2u)\|_{L^2(I_j;L^{\frac{2n}{n+2}}_x)}+
\||u|^2u\|_{L^{\frac{4}{3}}(I_j;L^{\frac{2n}{n+2}})}.
\end{equation}
Since $s_0\leq s<2$, by the fractional chain rule and H\"older's
inequality, we can control the nonlinearity as follows
\begin{equation}\label{5.20}
\big\||\nabla|^{s-1}(|u|^2u)\big\|_{L^2(I_j;L^{\frac{2n}{n+2}}_x)}\leq
\big\||\nabla|^{s-1}u\big\|_{L^{\infty}(I_j;L^{\frac{2n}{n-2}}_x)}
\|u\|^2_{L^{4}(I_j;L^{n}_x)},
\end{equation}
while by Sobolev embedding and interpolation
\begin{equation}\label{5.21}
\|u\|_{L^{4}(I_j;L^{n}_x)}\leq \|u\|_{L^{4}(I_j;\dot
H^{\sigma}_{p})}\leq\|u\|^{\theta}_{L^{4}(I_j;\dot
H^{s+\frac{1}2}_{\frac{2n}{n-1}})} \|u\|^{1-\theta}_{M(I_j)},
\end{equation}
where
\begin{equation*}
\begin{split}
\theta&=\frac{2n-5}{4s-1},\qquad \sigma-\frac np=-1,\\
\sigma&=\Big(s+\frac12\Big)\theta+\frac{5-n}4(1-\theta)=\frac{(2s+1)(2n-5)-(n-5)(2s-n+2)}{2(4s-1)}>0,\\
\frac1p&=\frac{n-1}{2n}\theta+\frac14(1-\theta)=\frac{(n-1)(2n-5)+n(2s-n+2)}{2n(4s-1)}.\\
\end{split}
\end{equation*}
On the other hand, we have
\begin{equation}\label{a5.19}
\||u|^2u\|_{L^{\frac{4}{3}}(I_j;L^{\frac{2n}{n+2}})}=\|u\|^3_{L^{4}(I_j;L^{\frac{6n}{n+2}})}\leq
\|u\|_{L^{4}(I_j;L^{n})}\|u\|^2_{L^{4}(I_j;L^{4})}.
\end{equation}
and
\begin{equation}\label{a5.20}
\|u\|_{L^4(I_j;L^4)}\lesssim\|u\|_{L_t^\infty(I_j;L^4)}^\frac{n-4}n\|u\|_{L_t^\frac{16}n(I_j;L^4)}^\frac4n\lesssim\|u\|_{S^s(I_j)}.
\end{equation}
This together with \eqref{5.19}-\eqref{5.21} yields that
\begin{equation*}
\|u\|_{S^s(I_j)}\lesssim\|\langle\nabla\rangle^su(t_j)\|_{L^2_x}+
\delta^{2(1-\theta)}\|u\|^{1+2\theta}_{S^s(I_j)}+
\delta^{1-\theta}\|u\|^{2+\theta}_{S^s(I_j)}.
\end{equation*}
A standard continuity argument yields that
\begin{equation*}
\|u\|_{S^s(I_j)}\lesssim\|\langle\nabla\rangle^su(t_j)\|_{L^2_x}\leq
C(\|u_0\|_{H^s}),
\end{equation*}
provided we choose $\delta$ sufficiently small depending on
$\|u_0\|_{H^s}$. Summing over all subintervals $I_j$, we have that
\begin{equation}\label{5.22}
\|u\|_{S^s(\R)}\lesssim C(\|u_0\|_{H^s}).
\end{equation}
To prove the asymptotic completeness, we need to prove that there
exist unique $u_\pm$ such that
\begin{equation*}
\lim_{t\rightarrow\pm\infty}\|u(t)-e^{it\Delta^2}u_\pm\|_{H^s_x}=0.
\end{equation*}
By time reversal symmetry, it suffices to prove the claim for
positive times only. For $t>0$, we will show that
$v(t):=e^{-it\Delta^2}u(t)$ converges in $H^s_x$ as
$t\rightarrow+\infty$, and $u_+$ to be the limit. In practice, we
can use Duhamel's formula to get
\begin{equation}\label{5.23}
v(t)=e^{-it\Delta^2}u(t)=u_0-i\int_0^te^{-i\tau\Delta^2}(|u|^2u)(\tau)\mathrm{d}\tau.
\end{equation}
Moreover, for $0<t_1<t_2$, we have
\begin{equation*}
v(t_2)-v(t_1)=-i\int_{t_1}^{t_2}e^{-i\tau\Delta^2}(|u|^2u)(\tau)\mathrm{d}\tau.
\end{equation*}
Using the Strichartz estimate, we derive that
\begin{equation*}
\begin{split}
\|v(t_2)-v(t_1)\|_{H^s_x(\R^n)}&=\Big\|\int_{t_1}^{t_2}e^{-i\tau\Delta^2}(|u|^2u)(\tau)\mathrm{d}\tau\Big\|_{H^s_x(\R^n)}\\&
\leq
\||\nabla|^{s-1}(|u|^2u)\|_{L^2([t_1,t_2];L^{\frac{2n}{n+2}}_x)}+
\||u|^2u\|_{L^{\frac{4}{3}}(I_j;L^{\frac{2n}{n+2}})}.
\end{split}
\end{equation*}
Arguing similarly as before, the above one can be controlled by
\begin{equation*}
\begin{split}
\|u\|^{2(1-\theta)}_{M([t_1,t_2])}\|u\|^{1+2\theta}_{S^s([t_1,t_2])}+\|u\|^{1-\theta}_{M([t_1,t_2])}\|u\|^{2+\theta}_{S^s([t_1,t_2])}.
\end{split}
\end{equation*}
Therefore, it follows from \eqref{5.1} and \eqref{5.22} that
\begin{equation}
\begin{split}
\|v(t_2)-v(t_1)\|_{H^s_x(\R^n)}\rightarrow0 \quad\text{as} \quad
t_1, t_2\rightarrow+\infty.
\end{split}
\end{equation}
As $t$ tends to $+\infty$, the limitation of \eqref{5.23} is well
defined. In particular, we find that
\begin{equation*}
\begin{split}
u_+=u_0-\int_0^\infty e^{-i\tau\Delta^2}(|u|^2u)(\tau)\mathrm{d}\tau
\end{split}
\end{equation*}
which is nothing but the asymptotic state. This concludes the proof
of Theorem \ref{the1}.
\end{proof}

{\bf Acknowledgements:}\quad The authors would like to thank Prof.
Jean-Claude Saut and  B. Pausader  for their invaluable comments
and suggestions. The authors were supported by
the NSF of China under grant No.11171033, 11231006.

\begin{center}

\end{center}
\end{document}